\documentclass [12pt]{article}
\input{artadj.sty}
\usepackage{latexsym}
\usepackage{graphicx}

\usepackage{amsfonts}
\usepackage{amsmath, amssymb}

\renewcommand{\thesection}{\arabic{section}}

\newcommand{\limn}{\lim_{n\rightarrow\infty}}
\newcommand{\limm}{\lim_{m\rightarrow\infty}}
\newcommand{\limk}{\lim_{k\rightarrow\infty}}

\newcommand{\limtau}{\lim_{\tau\rightarrow\infty}}

\newcommand{\goto}{\rightarrow}

\newcommand{\bo}{{\bf 1}}

\newcommand{\beqa}{\begin{eqnarray*}}
\newcommand{\eeqa}{\end{eqnarray*}}

\newcommand {\beq} {\begin{equation}}
\newcommand {\eeq} {\end{equation}}
\newcommand {\bear} {\begin{eqnarray}}
\newcommand {\eear} {\end{eqnarray}}
\newcommand {\bears} {\begin{eqnarray*}}
\newcommand {\eears} {\end{eqnarray*}}   

\newcommand {\done} {\quad\vrule height4pt WIDTH4PT}

\newcommand {\barr} {\begin{array}}
\newcommand {\earr} {\end{array}}

\newtheorem{thm}{Theorem}[section]
\newtheorem{lem}[thm]{Lemma}

\newtheorem{cor}[thm]{Corollary}

\graphicspath{%
    {converted_graphics/}% inserted by PCTeX
    {/}% inserted by PCTeX
}
\begin{document}
\bibliographystyle{plain}

\title{Classification of Discrete-Time Queues\\ 
}
%\date{}
%\maketitle
\author{ Muhammad El-Taha \\
                  Department of Mathematics and Statistics\\
                       University of Southern Maine\\      
                      96 Falmouth Street\\
                            Portland, ME  04104-9300\\
                     Email:el-taha@maine.edu
}
\date{ }

\maketitle

\ls{1}

\noindent{\bf Abstract.}
In this article we classify  discrete-time queues  based  on scheduling rules and observation epochs combinations. This  classification leads to {\em coherent}, {\em sub-coherent}, and {\em super-coherent}  systems when  {\em observed}  waiting times are, respectively equal to, less than, or larger than  {actual} waiting times. We then explore the consequences of this classification. Specifically, 
we discuss invariant properties of {\em coherent} systems including  queue-lengths, waiting times, servers' busy times, busy periods, Pollaczek-Khinchine formula,  and other common characteristics. An important consequence is that a performance characteristic of a system with specific scheduling rule and observation epoch combination extends to the entire class. 
An unresolved issue  in the literature is the  assertion that Little's law does not apply for  discrete-time queues that incorporate certain scheduling rules. Using this classification,  we reconcile the generality of Little's law and its applicability to all discrete-time queues regardless of  scheduling rules. 

%We start by 
%examining a general  statement and proof for {\em LL} that turns out to be helpful in
%resolving  the primary issue as to why {\em LL} does not appear to apply for certain discrete-time systems.

\bigskip

\bigskip

\noindent
%{\bf Mathematics Subject Classification:} Primary 60K25; Secondary 68M20; 90B36

\newpage
\ls{1.25}
%%%%%%%%%%%%%%%%%%%%%%%%%%%%%%%%%%%%%%%%%%%%%%%%%%%%%%%%%%%%%%%%%%%%%%%%%%%%%%%%%%%%%%%%%%%%%%%%%%%%%
\section{Introduction}\label{sec:intro}
%%%%%%%%%%%%%%%%%%%%%%%%%%%%%%%%%%%%%%%%%%%%%%%%%%%%%%%%%%%%%%%%%%%%%%%%%%%%%%%%%%%%%%%%%%%%%%%%%%%%%%%%%%%%%%%%%%%%%%%%%%%%%%%%%%%%%%%%%%%%%%%%

Consider a discrete-time queueing model where time $\tau=0,1,2,\ldots$ is discrete, and where time slots are of equal unit length. Arrivals and departures can occur at the same discrete time instants. We refer to this as the actual system.
To keep track of the system state and transitions between states one needs to order the arrivals and departures. This leads to several scheduling rules ({\em SR})  like the early arrival ({\em EAS}) and late arrival ({\em LAS}) systems. Moreover, to track the system state, (e.g., queue-length) one observes the system at slot edges or slot centers. Scheduling rules and observation epochs will be discussed in Section~\ref{sec:prel}. When considered together {\em SR} and observation epochs generate systems that are {\em coherent} and others that are {\em incoherent} in the sense that the customers' observed waiting times in these {\em incoherent} systems are  unequal to the actual system waiting times.

It is well-known that Little's law ({\em LL}) applies to discrete-time systems at great level generality. However, 
applying  discrete-time {\em LL} to queues with scheduling rules ({\em SR}) has not been fully explored. There are assertions in the literature that there are instances when {\em LL}
does not apply even when all assumptions are met. A primary motive in this article is to explore this issue and provide recommendations that remedy this situation.

Hunter~\cite{Hun83b} discusses three {\em SR}, the early arrival, late arrival with immediate access, and late arrival with delayed access. 
He shows that for {\em EAS} and LAS models the waiting time distribution function is the same for all three {\em SR} for G/Geom/1 and B/Geom/1 FCFS  models. But to derive his results, Hunter adopts the convention that (see page 228) for systems with {\em LAS-IA} he ``counts the number of service time positions spent in the system"; for {\em LAS-DA}  he ``counts the  completed number of time slots spent in the system"; and for the {\em EAS} system ``both methods give the same waiting time". It appears that  Hunter assumes apriori that the waiting time distribution function should be the same and the choice of multiple ways to count  service times is to achieve that outcome.
Desert and Daduna~\cite{Des02} page 85 conclude that ``{\em LL} can be directly applied to a discrete-time queueing system in the late arrivals case only". See also page 74 of their article. This is troubling as it is well known the {\em LL} applies universally to all discrete-time systems. 
Moreover, Desert and Daduna~\cite{Des02} study several {\em SR} ({\em EAS}, {\em LA-DF}, {\em LA-AF}) and conclude that waiting time distribution function is the same  for these {\em SR}. They also discuss  {\em LL} and conclude that (see pages 84-85) for certain observation  epochs {\em LL} does not hold in the sense that the {\em LL} cannot be applied to the {\em EAS} case. 
 Fiems and Bruneel~\cite{Fie02} consider {\em LL} for discrete-time systems as a special case of the continuous time version. They consider  a discrete time  equivalent  of the  {\em LL}  result in continuous time. They also study the effect of arrival/departure rescheduling  on {\em LL}. The discretized version leads to instances (see  equations (2) and (5) in  their paper) where it appears
 $L$ and $\lambda W$ are unequal.
Dattatreya and Singh~\cite{Dat05} discuss relationships among different discrete-time models that arise  in telecommunication networks. They study the {\em EAS} model and compare mean values (queue-lengths and response times) using slot edges and slot centers. They conclude (equation (17) in their paper) that {\em LL} applies when {\em L} is computed as slot centers, but not when {\em L} is computed at slot edges.
These articles are  examples of researchers' discomfort with applying {\em LL} and illustrations of the need for more clarity as to when and how {\em LL} applies to discrete-time queues.

Scheduling rules in discrete-time queues  have been addressed by several authors. Gravey and Hebuterne~\cite{Gra92} study the simultaneity of arrivals and departures in discrete-time queues,  the need for {\em SR}, and the effect of these {\em SR} on systems performance. 
Chaudhry et al.~\cite{Cha96} study  discrete-time models using {\em EAS} and {\em LAS-DA} at the random observer (slot edges) and outside observer (slot centers) epochs.
El-Taha et al.~\cite{Elt97a} introduce  the {\em LA-DF}  scheduling rule to study insensitivity of symmetric discrete-time queues.
Chaudhry~\cite{Cha00} gives a thorough discussion of the {\em EAS}, {\em LAS-DA} and {\em LAS-IA} scheduling rules.
Daduna~\cite{Dad01} discusses other issues that result from using scheduling rules. He mentions instances where {\em BASTA} does not hold even with Bernoulli arrivals in the sense  that Bernoulli arrivals should see time averages at the random observer epochs.
Articles that apply {\em LL} include Goswami~\cite{Gos14}, Bruneel et al.~\cite{Bru16}, Bruneel and Kim~\cite{Bru93} and Chaudhry and Gupta~\cite{Cha97}, among others.

One can identify  five popular  scheduling rules and six  observation epochs (Section~\ref{sec:prel}). This creates a large number of potential  instances that may require attention. This article creates a new framework that addresses these issues and provides a resolution. Specifically, the  contributions of this article are:

1. Classify discrete-time queues based on  actual and observed waiting times. This is a novel idea that has not been utilized before in the literature. We examine all {\em SR} and observation epochs arrangements and classify them into {\em coherent}, {\em sub-coherent} or {\em super-coherent} systems based on  whether a {\em SR}/observation epoch combination leads to the observed  waiting time being equal, less than or larger than the actual waiting. We shall carefully define  what we mean by {\em coherent}, {\em sub-coherent}, {\em super-coherent}, {\em actual}, and {\em observed}  in Sections~\ref{sec:prel} and \ref{sec:class} of this  article.
 
2. Show that members of each class  share invariant characteristics, thus can be studied together as one unit.  In particular, we show that all {\em coherent} systems  share the same stationary distribution, the same waiting time distribution function, the same  servers' utilization factors, the same busy periods, and many more. Thus, studying one {\em coherent} system will give results that are shared by the entire class.
Contrast this with the literature,  where each {\em SR}/observation epoch  arrangement is studied independently.

3. We address an unresolved issue that {\em LL} does not appear to hold for some discrete systems. In particular, we show  that when {\em LL} does not appear to hold, it is because, for incoherent systems, either  $L$ or  $W$ is computed using actual  values while the other is computed using observed values. Specifically, if we define $W$ using the actual (observed), we need to consistently, define $L$ using the actual (observed) values. We also give relationships between $L$ and $W$ when one is based on  observed values and the other is based on actual values.

The rest of the article is organized as follows. 
In Section~\ref{sec:prel} we give preliminary discussion that includes {\em SR}, observation epochs and  motivation where we discuss issues with applying {\em LL} to queues with {\em SR}. We also discuss waiting times in the presence of {\em SR} and perform two experiments and observe a relation between actual and observed waiting times   for some combinations of {\em SR} and observations epochs.
In Section~\ref{sec:class}  we  formalize our observations in Section~\ref{sec:prel} and classify discrete-time queues with {\em SR} into {\em coherent} and {\em incoherent} systems and discuss each type of system. 
Moreover, we prove that the queue-length distribution is invariant for all {\em SR} within {\em coherent} systems.
In Section~\ref{sec:LL},  we review {\em LL} for general discrete systems where no ordering of arrivals and departures is assumed. The proof  is deferred to the appendix. Then, we show that {\em LL} applies to all systems when  waiting times are  appropriately defined. We then apply {\em LL} to {\em coherent} and {\em incoherent} systems.
In Section~\ref{sec:LLbdq}  we consider the special case of birth death processes with focus on {\em coherent} and {\em sub-coherent} systems.
In Section~\ref{sec:invar}  we further characterize {\em coherent} systems by exploring their  invariant properties. We focus on busy servers and busy periods.
In Section~\ref{sec:HlG}  we apply  $H=\lambda G$, an extension of {\em LL},  to give  a proof for Pollaczek-Kinchine formula for all {\em coherent} $B/G/1$ models.
In Section~\ref{sec:cr}  we give concluding remarks.

\section{Preliminaries: Queues with Scheduling Rules}\label{sec:prel}

In this section, we discuss preliminaries that are needed in Section~\ref{sec:class}. Specifically, we discuss what we mean by the actual/standard system, introduce scheduling rules and observation epochs, and discuss our motivation and known issues that arise with {\em LL} when applied to queues with {\em SR}.  We also give examples that compare actual vs observed waiting times.

\subsection{The Actual/Standard System}

We define a discrete-time {\em actual} or {\em standard}  system as one where arrivals and departures take place at discrete-times $\tau$ (slot edges) exactly. This means that  more than one event type can take place at the same time instant. Sometimes this is all we need to obtain certain results as we do with {\em LL} in Section~\ref{sec:LL}. However, quite often, we need to order the events (arrivals and departures) such as when we need to track state transition probabilities. This leads to a variety of {\em SR} that are described next.

\subsection{Scheduling Rules}\label{sec:sr}

Here we describe  discrete-time queueing models where time is divided, w.l.o.g., into slots of unit length.  We assume a discrete-time scale $\{0,1,\ldots, \tau,\ldots\}$ so that a time slot is an interval of the time 
$(\tau, \tau+1]$.
The system is driven by two event types: arrivals and departures (service completions). We assume that the events, i.e., arrivals  and or departures   occur at the boundary of the time slots. The order of scheduling arrivals and departures in a time slot (equivalently at a time instant) leads to different {\em SR} that we describe below.
 Depending on the behavior of the physical system the order of potential arrivals and departures at any given slot vary significantly. This results in  various {\em SR}, also referred to as scheduling systems or waiting room management policy. Let $A$ and $D$ denote a potential  arrival or departure, respectively, at time $\tau$. We follow the notation setup as in Hunter~\cite{Hun83b}, Chaudhry et al~\cite{Cha96} and Desert and Daduna~\cite{Des02}. 
%\bigskip

In the early arrival system ({\em EAS}),  potential arrivals in a time slot are scheduled to occur before potential departures.
Specifically, a potential arrival at time $\tau$  occurs in $(\tau,\tau+)$, and a potential departure at time $\tau$ occurs  in $(\tau-,\tau)$. 
That is, $\tau-<D<\tau<A<\tau+$. 
Moreover, if an arrival finds an idle server, it goes into service immediately  and can potentially depart in the same time slot.
In the late arrivals system ({\em LAS}) the order of potential arrivals and departures is reversed so that  potential departures occur early in a time slot and potential arrivals occur at the end of the slot. 
 More specifically, a potential departure at time $\tau$  occurs in $(\tau,\tau+)$, and a potential arrival at time  $\tau$ occurs  in $(\tau-,\tau)$.
That is  $\tau-<A<\tau<D<\tau+$. 
Moreover, if an arrival at $\tau-$ finds an idle server and starts service immediately, it can potentially depart at the beginning of the next time slot at $\tau+$, the system is called immediate access ({\em IA}).
The late arrival with delayed access ({\em LAS-DA}) scheduling system is similar to the {\em LAS-IA} except that 
 an  arrival at $\tau-$ waits until the next  slot to start service at $\tau+1$, then the system is called delayed access ({\em DA}). 
 In the late arrivals departures-first ({\em LA-DF}) system both potential arrivals and departures occur late in the slot, so that $\tau--<D<\tau-<A<\tau$. An arrival that finds an idle server stars service at $\tau$.
%\bigskip
%
 Finally, in the late arrivals arrivals-first ({\em LA-AF}) system both potential arrivals and departures occur late in the slot, so that $\tau--<A<\tau-<D<\tau$.  An arrival that finds an idle server stars service at $\tau$. These scheduling rules are depicted in Figure 1.

\begin{figure}[tbp] % float placement: (h)ere, page (t)op, page (b)ottom, other (p)age
  \centering
  % file name: C:/eltaha25/papers/Discrete-dtspa/A-Classification-LL-InvarianceDiscreteTime/August25/SRfig.eps
  \includegraphics[bb=0 0 1376 769,width=5.67in,height=3.17in,keepaspectratio]{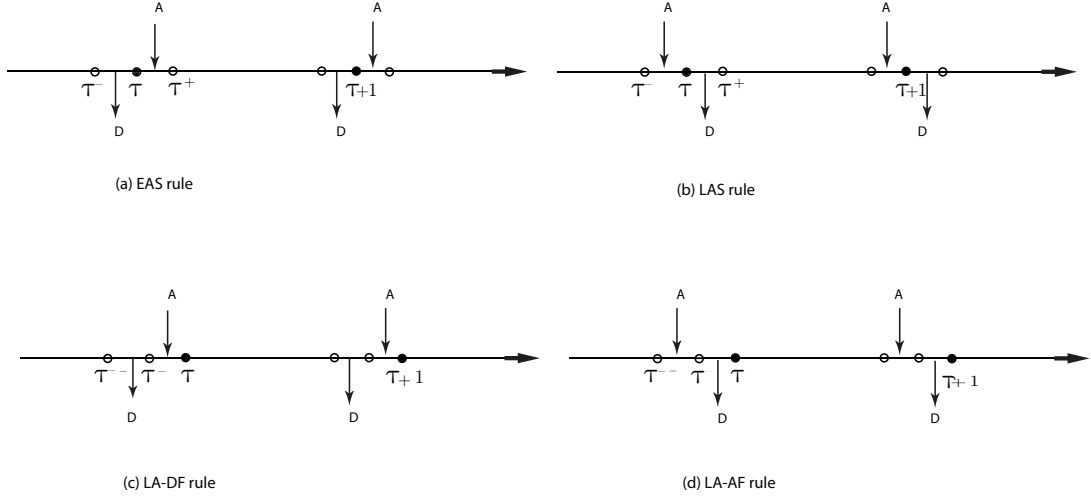}
  \caption{Scheduling rules representation}
  \label{fig:SRfig}
\end{figure}

For details about {\em EAS}, {\em LAS-IA} and {\em LAS-DA}  scheduling regimes one may consult Hunter~\cite{Hun83b} and  Chaudhry~\cite{Cha00}. The {\em LA-DF} is first introduced by El-Taha et. al. ~\cite{Elt97a}. Moreover, the {\em LA-DF} and {\em LA-AF} are discussed by Daduna~\cite{Dad01}, and Desert and Daduna~\cite{Des02}.
%In this article, we follow the notation setup as in Hunter~\cite{Hun83b}, Chaudhry et al~\cite{Cha96} and Desert and Daduna~\cite{Des02}. 
For more discussions about  these {\em SR} one can consult El-Taha~\cite{Elt23} and Gravey and Hebuterne~\cite{Gra92} for a reference on this. 

\subsection{Observation Epochs}

An important feature of discrete-time systems is that we can obtain time average measures using more than one observation epoch.  Contrast this with continuous-time systems where time average measures are obtained by averaging continuously over time. 
Here we describe six epochs. 
The most natural time-average measure is to average over  slot edges or random observer epochs. Another common  epoch is the outside observer  or slot center epochs.  It is also common to use scheduled potential pre-arrival, post-arrival, pre-departure, and post departure epochs. Note that every value of $\tau$ is a potential (vs actual) arrival and departure epoch.
Let $u(\tau)$ be an observation epoch. We have six observation epochs for each of the five {\em SR}. In Table~\ref{table-epoch} we list those values.

\begin{center}
\noindent
\begin{tabular}{|c|cccccc|}
\multicolumn{6}{c}{\bf Table~\ref{table-epoch} Possible values of $u(\tau)$} \\ 
\multicolumn{6}{c}{} \\
 \hline 
 & Random & Outside &Poten. Pre& Poten. Post&Poten. Pre&Poten. Post  \\
& Observer& Observer&Arrival&Arrival&Departure&Departure\\       
\hline
$EAS$    &$\tau$ &$\tau-.5$ &$\tau$ &$\tau^+$&$\tau^-$&$\tau$ \\
LAS-IA    &$\tau$ &$\tau-.5$ &$\tau^-$ &$\tau$&$\tau$&$\tau^+$ \\
LAS-DA  &$\tau$ &$\tau-.5$ &$\tau^-$ &$\tau$&$\tau$&$\tau^+$ \\
 LA-AF   &$\tau$ &$\tau-.5$ &$\tau^{--}$ &$\tau^-$&$\tau^{-}$&$\tau$ \\
 LA-DF   &$\tau$ &$\tau-.5$ &$\tau^{-}$ &$\tau$&$\tau^{--}$&$\tau^-$ \\
\hline
\end{tabular}
\label{table-epoch}
\end{center}
%\vspace{5mm}
%
Note that, for each scheduling rule, the system state is observed at one of these epochs and then system's characteristics are obtained. One issue  is how to measure waiting time in the system. It has been suggested that  time in system can be measured using service positions and/or service slots, see Hunter~\cite{Hun83b}, without regard to the observation epoch. There is an issue with counting service slots since we only observe the system at discrete-time points. However, counting service slots is equivalent  to counting service positions at  slot centers. This works because arrivals and departures take place at slot edges. We believe counting service positions at slot centers makes more sense, however, since  slot centers are possible observation epochs. In this article we shall use service positions at slot centers instead of service slots.

\subsection{Motivation}

In this subsection we point out  instances in the  literature where it is stated that  {\em LL} does not appear to hold.
Consider a stable B/Geom/1 model. Arrivals follow a Bernoulli process such that the probability of an arrival in any slot is $0<\alpha <1$. Service times  are geometric with  parameter 
$0<\beta<1$, so that 
$
P(S=n)=(1-\beta)^{n-1} \beta \;\; (n= 1,2,\ldots )\;,
$
and $E[S]=1/\beta\;$. Let $\rho=\alpha/\beta <1$. 
Consider system characteristics at the slot edges (random observer) and slot centers (outside observer).
The superscripts $R$ and $C$ will be used to represent random observer (slot edges) and outside observer (slot centers) respectively.
Note that for this system, and regardless of scheduling rule, it is well established (Hunter~\cite{Hun83b}) that the mean  waiting time in the system is given by $W=\frac{1-\alpha}{\beta-\alpha}.$   

Now, as noted in  Desert and Daduna~\cite{Des02}, apply {\em LL} with $\lambda=\alpha$ to {\em EAS} and {\em LAS-IA} systems. Consider the {\em EAS} at the random observer epochs, then it is well-known (e.g., El-Taha~\cite{Elt23}) that $L^R =\frac{\alpha(1-\beta)}{\beta-\alpha}$. But $L^R \not=\lambda W$. Similarly, consider the {\em LAS-IA} system observed at the slot centers. Then 
for this system $L^C =\frac{\alpha(1-\beta)}{\beta-\alpha}$. However, $L^C \not=\lambda W$.% These observations are made by Desert and Daduna~\cite{Des02}.
 The conclusion is that {\em LL} does not appear to apply for {\em EAS} system at the random observer epochs and {\em LAS-IA} at the outside observer epochs. On the face of it this appears, as noted in the literature, to contradict $LL$.

Some combinations of {\em SR}/observation epochs have an unexpected effect on the application 
of {\em LL}. In the next subsection we examine three examples to shed light on issues with evaluating waiting times in queues when we factor in {\em SR} and observation epochs.

\subsection{Waiting Times in  Queues with Scheduling Rules}\label{subsec:ex}

 We continue our motivation and consider three examples that illustrate the issue with waiting times in discrete-time queues. 
We assume that in these queueing models service times are greater or equal to one, i.e., no zero service times.
Let  $A_k, D_k$ be the actual system's $k^{th}$ arrival and departure epochs, so that
\begin{equation}\label{eq:actualW}
W_k=\sum_{\tau=0}^\infty {\bf 1}\{A_k < \tau \leq D_k\}=D_k-A_k\;.
\end{equation}

 Now, let $A'_k,D'_k$ be the scheduled $k^{th}$ arrival and departure epochs associated with selected  {\em SR}.  For a given {\em SR}, arrival and departure instants are related to the actual arrival and departure instants as follows: 
\begin{equation*}\label{eq:A'}
A'_k=\left\{\begin{array}{lll}
        A_k^- &\mbox{for LAS-IA, LAS-DA, LA-DF}\;;\\\\
        A_k^{--} & \mbox{for LA-AF}\;;\\\\
	A_k^+ & \mbox{for EAS}\;.
\end{array}  \right.              
\end{equation*}
   Similarly,
\begin{equation*}\label{eq:D'}
D'_k=\left\{\begin{array}{lll}
        D_k^- &\mbox{for EAS, LA-AF}\;;\\\\
        D_k^{--} & \mbox{for LA-DF}\;;\\\\	
	D_k^+ & \mbox{for LAS-DA}\;;\\\\
	(D_k-1)^{+} & \mbox{for LAS-IA}\;.
\end{array}  \right.              
\end{equation*}
Let the observed waiting time  at slot edges of $k^{th}$ arrival be  $W_k^o=\sum_{\tau=0}^\infty {\bf 1}\{A'_k < \tau \leq D'_k\}$.

\noindent
{\bf Example 1.}
Consider a single server queue where, say   $k^{th}$ customer, with one unit of service $S_k=1$ arrives at an idle server at some time point $\tau$ so that $W_k=1$. Let us consider the random observer (slot edges) epochs. Consider the {\em LAS-IA} model where, $A'_k=\tau^-$ and $D'_k=\tau^+$ . Here we see that
 $D'_k-A'_k=0\not=W_k=1$.  
However, the observed waiting time  $W_k^o=1= W_k$ giving the correct time in the system. Now, consider the {\em EAS} model where $A'_k=\tau^+$ and $D'_k=(\tau+1)^-$ and note that  $D'_k-A'_k=1=W_k$,   giving the correct answer. However, the observed waiting time  $W^o_k=0\not= W_k=1$ giving the incorrect time in the system. Using  {\em LAS-DA} system we see that $D'_k-A'_k=1=W_k$, and $W^o_k=2$.
 Both  definitions  $D'_k-A'_k$ and  $W^o_k$,  give the  correct $W_k=1$ for both {\em LA-DF} and {\em LA-AF} rules. 
\bigskip

\noindent
{\bf Example 2.}
Consider a $G/G/\infty$ model. In this model all customers arrive to find an idle server and the waiting time in the system is exactly the service time. Let us consider the first arrival and assume this customer arrives at time $\tau=1$ and requires $S_1$ units of service. Now, consider the five {\em SR} using service positions at slot edges and  service slots counted at slot centers. 

 Counting service positions (equivalently, observing the system at slot edges), we obtain the correct value $W_1=W^o_1=S_1$ for all {\em SR} except {\em EAS} and {\em LAS-DA} 
where we obtain the incorrect waiting time  where $W^o_1=S_1-1$ and $W^o_1=S_1+1$, respectively. On the other hand, observing the system at slot  centers (outside observer),  we obtain the correct value $W_1=W^o_1=S_1$ for all {\em SR} except {\em LAS-IA} where we obtain the incorrect waiting time where  $W^o_1=S_1-1$. 

Depending on the  {\em SR} and observation epochs combination we end up with $W^o=E[S]$, $W^o=E[S]-1$ or 
$W^o=E[S]+1$. This incorrect evaluation of $W$  leads to incorrect application of {\em LL}, when $L$ is based on the actual arrival and departure  times. 
\bigskip

\noindent
{\bf Example 3.} In this example we show how incorrectly counting waiting times extends to busy periods.
Consider a single server model with two arrivals such that $A_i=1,3$ and $S_i=5,4$ for $i=1,2$. The process repeats every $10$ units. We focus on the first two arrivals and their service times in the system which is the first busy period $B_1=9$  Observing the system at slot edges (random observer), we obtain the correct busy period where  $B^o_1=9$ for all {\em SR} except {\em EAS} and {\em LAS-DA} where we obtain the incorrect busy period with $B^o_1=8$ and $B^o_1=10$, respectively.
On the other hand, counting service at slot centers (outside observer), we obtain the correct busy period with $B^o_1=9$ for all {\em SR} except {\em LAS-IA} where we obtain the incorrect busy period with $B^o_1=8$. 
\bigskip

The observations in these examples are summarized in Table~\ref{table-count}.  A yes in the table is for combinations of {\em SR}/observation epochs that give correct waiting times.

\begin{center}
\noindent
\begin{tabular}{|c|ccc|}
\multicolumn{4}{c}{\bf Table~\ref{table-count} Counting Actual Waiting Times} \\ 
\multicolumn{4}{c}{} \\
 \hline 
 & Counting Service  & Counting Service &   \\
& Positions at Slot Edges & Positions at Slot Centers & \\     
\hline   
$EAS$       &no &yes & \\
LAS-IA      &yes  &no& \\
LAS-DA   &no &yes &\\
 LA-AF  & yes&yes& \\
 LA-DF  & yes&yes& \\
\hline
\end{tabular}
\label{table-count}
\end{center}
\vspace{5mm}

In {\em actual} discrete-time systems, where arrivals and departures occur exactly at  the observed integer
values $\tau$, we see that the middle part of  (\ref{eq:actualW}) counts the number of service positions that the customer is in the system while r.h.s. of (\ref{eq:actualW}) counts the number of slot centers the customer spends in the system. More importantly,  the time in the system is the same in both cases.  With scheduling rules, observed and actual waiting times are not necessarily the same. This is addressed in the next section. 

\section{Classification of  Discrete-Time Queues}\label{sec:class}

In this section we classify discrete-time systems  with respect to {\em SR} and observation epochs combinations.
For our purposes, a discrete-time queueing system has three basic features. The first feature is the {\em SR} which can be any of {\em EAS}, {\em LAS-IA}, {\em LAS-DA}, {\em LA-AF}, {\em LA-DF}. One can think of others, but these are the ones considered in the literature. The second feature is the state observation epoch needed to obtain the time-average system characteristics. In continuous time one observes the system continuously over time. In discrete-time queues with {\em SR}, we have  the random observer epochs (slot edges), outside observer epochs (slot centers), and  potential pre-arrival, post arrival, pre-departure, and post departure epochs.  The third feature is how to measure time in the system including  time in queue and time in service. In the literature, regardless of the observation epochs, there are two methods to measure time in system by either   counting time positions (e.g., slot edges or slot centers), or counting service slots. More importantly,  the service position/slot is not always the selected observation epoch, contributing to the issues mentioned earlier. By contrast we count time in the system at the observation epochs.   

\subsection{Assumptions}

With {\em SR} that schedule arrivals and departures around slot boundaries, observed and actual waiting times are not always equivalent. To resolve this issue, we introduce the concept of {\em coherent} systems.
Consider any discrete-time queueing system with possible multiple {\em SR} and observation epochs. Let $u(\tau)\in \{\tau,\tau-.5,\tau^{--},\tau^-,\tau^+,\tau^{++}\}$ be an observation epoch, and $\{(A_k,D_k), k \geq 1\}$be the input data.
 At the heart of it, it turns out that
the actual customer waiting 
 $ W_k=D_k-A_k$ and the observed waiting time (waiting time as observed at the observation epochs) are not equal for all {\em SR}/observation epochs combinations.
This observation motivates our next definition.
\bigskip
  
\noindent
{\bf Definition 1.}   
A discrete-time queueing  system is said to be {\em coherent} if for all customers, $k=1,2,\ldots $,
\begin{equation}\label{defeq:coherent2}
 \sum_{\tau=1}^{\infty}\bo\{A'_k < u(\tau)\leq D'_k\}=\sum_{\tau=1}^{\infty}\bo\{A_k < \tau \leq D_k\}\;.
\end{equation}
Otherwise, the system is said to be {\em incoherent}. \done
\bigskip 

This definition says that a discrete-time queueing system is {\em coherent} if, for every customer, the observed waiting time is equal to the actual waiting time. In other words, when you measure how long a customer has been in the system based on a specific scheduling rule and observation epoch, the result matches the actual time elapsed between their arrival and departure. An incoherent system is one where the observed waiting time does not equal the actual waiting time.
\pagebreak

\noindent
{\bf Remarks.} 

(i) Let
\[
W^o_k= \sum_{\tau=1}^{\infty}\bo\{A'_k < u(\tau)\leq D'_k\}
\]
be the  {\em observed}  waiting time associated with some {\em SR}/observation epochs combinations. Then (\ref{defeq:coherent2}) says that a system is  {\em coherent} if for all arrivals $k=1,2,\ldots $,
\[
W^o_k=W_k\;.
\]
Otherwise, the system is said to be {\em incoherent}.

(ii) Note that for the {\em actual} system for all $k=1,\ldots $
\[
\sum_{\tau=1}^{\infty}\bo\{A_k < \tau \leq D_k\}=\sum_{\tau=A_k}^{D_k}\bo\{A_k < \tau \leq D_k\} =D_k- A_k:=W_k
\]

 (iii) For systems with {\em SR} and random observation epochs, i.e., $u(\tau)=\tau$ we have
\[
\sum_{\tau=1}^{\infty}\bo\{A'_k < \tau \leq D'_k\}=\sum_{\tau=A_k}^{D_k}\bo\{A'_k < \tau \leq D'_k\} 
\]
This sum equals  the actual $W_k$ of corresponding actual system for only {\em coherent} systems.

\begin{lem}\label{lem:wt-coh} The  waiting time distribution function is invariant  with respect to {\em SR} for all {\em coherent} systems.
\end{lem}
{\bf Proof.}
The proof follows from observing that for any sample path,  all $\{W_k, k\geq 1\}$ for any  {\em coherent} system are the  same as that of the actual system.\done

Incoherent systems have been recognized and studied in the literature, see for example Hunter~\cite{Hun83b}.
The examples in Subsection~\ref{subsec:ex} suggest that there are two types of {\em incoherent} systems. One is where the observed waiting time is shorter the actual waiting time, and the other is where the observed waiting time is longer the actual waiting time.

\noindent
{\bf Definition 2.} Consider any discrete-time queueing system with possible multiple {\em SR} and observation epochs. Then the system is said to be {\em sub-coherent} if for some arrivals $k=1,2,\ldots $, 
\[
\sum_{\tau=1}^{\infty}\bo\{A'_k< u(\tau) \leq D'_k\}< D_k-A_k:=W_k\;,
\]
i.e., if for some $k$,
\[
W^o_k<W_k\;,
\]
 the observed are shorter than the actual  waiting times. The system is said to be  {\em super-coherent} if for some arrivals $k=1,2,\ldots $,
\[
\sum_{\tau=1}^{\infty}\bo\{A'_k< u(\tau) \leq D'_k\}> D_k-A_k:=W_k\;,
\]
i.e.,  if  for some $k$,
\[
W^o_k>W_k\;,
\]
 the observed  are longer than the actual waiting times.
\done

Definition 2 says that   incoherent systems can be further classified into 
sub-coherent where  the observed waiting time is shorter than the actual waiting time, and 
super-coherent where the observed waiting time is longer than the actual waiting time.

 The following result holds on the sample paths.
\begin{lem}\label{pro:wo} For all $k=1,\ldots$,

 (i) if a queue is {\em sub-coherent} then
\[
W^o_k :=\sum_{\tau=1}^{\infty}\bo\{A'_k< u(\tau) \leq D'_k\}=W_k-1\;, and 
\]
(ii) if a queue is {\em super-coherent} then
\[
W^o_k := \sum_{\tau=1}^{\infty}\bo\{A'_k< u(\tau) \leq D'_k\}=W_k+1\;.
\]
\end{lem}
An immediate consequence of Lemma~\ref{pro:wo} is that Lemma~\ref{lem:wt-coh} applies to sub-coherent and super-coherent classes.
The following table identifies  {\em SR}/epochs combinations that lead to {\em coherent}, {\em sub-coherent} and {\em super-coherent}  systems.

\begin{center}
\noindent
\begin{tabular}{|c|cccccc|}
\multicolumn{6}{c}{\bf Table~\ref{table-coh} Coherent,  sub-coherent  and super-coherent  systems} \\ 
\multicolumn{6}{c}{} \\
 \hline 
 & Random & Outside &Poten. Pre& Poten. Post&Poten. Pre&Poten. Post  \\
& Observer& Observer&Arrival&Arrival&Departure&Departure\\       
\hline   
$EAS$       &sub &coh &sub &coh&coh&sub \\
LAS-IA      &coh  &sub& sub&coh&coh&sub\\
LAS-DA   &super &coh &coh &super&super&coh\\
 LA-AF  & coh&coh& coh&super&super&coh\\
 LA-DF  & coh&coh& sub&coh&coh&sub\\
\hline
\end{tabular}
\label{table-coh}
A `coh' indicates a coherent, `sub' indicates sub-coherent, and `super' indicates super-coherent systems.
\end{center}
Because (actual) pre and post event times are  subsets of the corresponding  potential pre and post event times, we have the following immediate result.
\begin{cor}
If a system is {\em coherent} at potential pre-arrival, pre-departure, post-arrival, post-departure  epochs, then it is {\em coherent} at the corresponding pre-arrival, pre-departure, post-arrival, post-departure  epochs.
\end{cor}

This result is useful in many regards, for instance  $BASTA$ relates pre-arrival customer average distribution function to time-average distribution functions. This Corollary says all {\em coherent}  systems have the same pre-arrival  customer-average  distribution function.

\subsection{Stationary Queue-Length Distribution of Coherent Queues}

In this subsection we show that the stationary distribution of the number of customers in the system is invariant  with respect to all {\em coherent} queues. 
The definition of {\em coherent} systems leads to the following result  about service times. 
\begin{lem} Let $S_k$ be the service requirement of the $k^{th}$ arriving customer. Then, for all {\em coherent} work conserving non-preemptive systems the actual and  observed service times, are equal. That is
 for all $k=1,2,\ldots$, $S_k=\sum_{\tau=1}^{\infty}\bo\{A''_k < u(\tau)\leq D'_k\}\;,$ where 
$A''_k$ is the start of service of the $k^{th}$ arrival.
\end{lem}
{\bf Proof.} For customers that wait, the service requirement and the time   in service are equal. For customers that find an idle server   $A_k''=A'_k$, so that the observed service time is the same as the  observed waiting time. The result follows from the definition of {\em coherent} systems.\done

This lemma implies that, for {\em coherent} systems, the total time a server is busy serving a customer  equals its service requirement regardless of how we count a unit of service. The consequence of this is that for {\em coherent} systems we have $W=W_q+ES$, where $W_q$ and $ES$ are, respectively, the mean waiting time in the queue and mean service time per customer. The same cannot be said about {\em incoherent} systems. Moreover, this result is useful is asserting that busy periods  are invariant for {\em coherent} systems (see Lemma~\ref{lem:busyP}). This is needed in proving the next result.

We know that for $\tau=1,\ldots $,   $L(u(\tau))$, is not invariant even for {\em coherent} systems.  On the other hand, we know that that the mean number of customers,  $L$, is invariant  for {\em coherent} systems. Here we show that the stationary distribution function is  invariant with respect to {\em coherent} systems.

For $n=0,1,\ldots$; $k=1,\ldots$, let $U_k$ be the start of $k^{th}$ busy cycle,
 $Y'(n,\tau)=\sum_{j=1}^{\tau}\bo\{L(u(j))=n\}$ $\left(Y(n,\tau)=\sum_{j=1}^{\tau}\bo\{L(j)=n\}\right)$ 
be  the  observed (actual) total time  in state $n$ during $(0,\tau]$;
 $C'_k(n)=Y'(n,U_{k+1})-Y'(n,U_k)$ $\left(C_k(n)=Y(n,U_{k+1})-Y(n,U_k)\right)$
 be the observed (actual) total time in state $n$ during $k^{th}$  busy cycle. Moreover, for all $n=0,1,\ldots,$ define the following limits when they exist:
\begin{eqnarray*}
\pi'(n)&=&\limtau\frac{Y'(n,\tau)}{\tau}\;;\\
\pi(n)&=&\limtau\frac{Y(n,\tau)}{\tau}\;;\\
C'(n)&=&\limm\sum_{k=1}^mC'_k(n)/m\;;\\
C(n)&=&\limm\sum_{k=1}^mC_k(n)/m\;.
\end{eqnarray*}
We interpret $\pi'(n)(\pi(n))$ as the observed (actual)  state-$n$ long-run frequency, and $C'(n)(C(n))$ as the observed (actual) long-run average number of visits to state $n$ per cycle. We assume that $\pi(n)$ is a proper distribution function.  

\begin{thm}\label{thm:df} The  queue-length stationary  distribution function $\pi'(.)$ is invariant with respect to all {\em coherent} systems. Specifically, $\{\pi'(n)\}$ and $\{\pi(n)\}$$\;, n=0,1,\ldots\;$ are identical.
\end{thm}
{\bf Proof.} 
First, we  
 show that  for all $n,$
\[
\pi'(n)=C'(n)/C'\;,
\]
where $C'$ is the busy cycle associated with the {\em coherent} class.  
Now, we apply a discrete-time counterpart of  $Y=\lambda X$ (See El-Taha and Stidham~\cite{Elt99}). Let $Y(\tau)=Y'(n,u(\tau))$, $X_k=C'_k(n)=Y'(n,U_{k+1})-Y'(n,U_k)$, so that $Y=\pi'(n)$, and $X=C'(n)$. Moreover, 
$\lambda=\limk\frac{k}{U_k}=\limtau \frac{A'(0,u(\tau))}{\tau}$, where $A'(0,u(\tau))$ counts the number of arrivals that see the system in state $0$, i.e., it counts the number of busy cycle starts  up to $u(\tau)$.  Now
\begin{eqnarray*}
C'(n)&=& X=\lambda^{-1}Y\;,\\
    &=&\limtau \frac{\tau}{A'(0,u(\tau))}\frac{Y'(n,\tau)}{\tau}\;\\
     &=&\limtau \frac{Y'(n,\tau)}{A'(0,u(\tau))}\;.\\
\end{eqnarray*}
Noting that $|u(\tau)-\tau| < 1$, $C'=\sum_{n=0}^{\infty} C'(n)=\limtau\frac{u(\tau)}{A'(0,u(\tau))}$. Now,
\begin{eqnarray*}
\pi'(n)&=&\limtau\frac{Y'(n,\tau)}{\tau}\;\\
        &=&\limtau \frac{Y'(n,\tau)}{A'(0,u(\tau))}\frac{A'(0,u(\tau))}{\tau}\;\\
     &=&C'(n)/C'\;.
\end{eqnarray*}
Similarly,  we can show that,
\[
\pi(n)=C(n)/C\;.
\]
It follows from Lemma~\ref{lem:busyP} that $C$, the mean cycle length, is invariant with respect to  {\em coherent} systems.
Therefore, $C=C'$ for all {\em coherent} systems. It remains to show that  $C'(n)=C(n)$ for all $n$.
Now, observe that for all {\em coherent} systems we have, 
\begin{equation}\label{eq:ql}
\sum_{j=U_{k}+1}^{U_{k+1}} \bo\{L(u(j)=n\}= \sum_{j=U_{k}+1}^{U_{k+1}} \bo\{L(j)=n\}\;, n=0,1,\ldots \;,
\end{equation}
i.e., for {\em coherent}, systems the number of visits of the queue-length process to state $n$ during a busy cycle is equal to the visits to state $n$ of the actual system. Therefore, by (\ref{eq:ql}) we have   $C'_k(n)=C_k(n)$ which implies  $C'(n)=C(n)$. This completes the proof. \done
  
In the literature {\em coherent} systems are studied individually, for example, Chaudhry et al.~\cite{Cha96}  study {\em coherent} versions  of the  {\em EAS} and {\em LAD-DA}  queues separately and compare their performance characteristics. Theorem~\ref{thm:df} is an important result in that it says  that all stable {\em coherent} systems share the same stationary distribution function.
Additionally, consider a recent article by Grassmann and Tavakoli~\cite{Gra19} where they study the queue-length distribution function for discrete-time {\em GI/G/1} queue using direct and indirect methods.
The queue-length is observed before any event (potential arrival or departure) takes place at time $\tau$. This is equivalent to the  outside observer (slot center) epochs. So, their results apply to all scheduling rules at slot centers except  {\em LAS-IA} which is not a {\em coherent} system at that epoch. Moreover,  for {\em LAS-IA} the random observer epoch falls after the potential arrival, but before the potential departure, so it is explicitly excluded from their study. Applying our results in this section, we see that their results should apply to all {\em coherent} systems including the {\em LAS-IA} at the random observer (slot edges) epoch.

In evaluating $L$ and $W$ we need consistency in counting the number of customers in the system and the time spent in the system for all arriving customers. 
This raises the question of how to correctly compute  $W$, and how to apply {\em LL} properly. This is addressed in the next section.
%\bigskip

%%%%%%%%%%%%%%%%%%%%%%%%%
\section{Little's Law for Discrete Systems}\label{sec:LL}

 To pinpoint the issues surrounding the application of {\em LL} to discrete-time queues we need a general rigorous proof of {\em LL}.
However, to the best of our knowledge, there is no proof of  discrete-time {\em LL} at the level of generality of the proof given by Stidham~\cite{Sti74} and El-Taha and Stidham~\cite{Elt99} for the continuous time {\em LL}.  Moreover,  the definitions of $L$ and $W$ play a role in how {\em LL} functions when {\em SR} are invoked.  In this section we review the discrete-time {\em LL} and give  a sample-path proof at a level of generality comparable to Stidham~\cite{Sti74}. It turns out that the  proof of {\em LL} for discrete-time systems is similar to that of the continuous time counterpart.

Consider a discrete-time scale $\tau=0,1,\ldots$. The input  
data are $\{(A_k,D_k), k \geq 1\}$, where $0 \leq A_k \leq A_{k+1} < 
\infty$, $A_k \leq D_k < \infty$, $k \geq 1$, and $A_k$ and $D_k$ are 
interpreted as the arrival
time and the departure time, respectively, of customer $k$. We assume $\{A_k, k\geq 1\}$  is a deterministic point process, i.e., a sample path of a stochastic process.
% with $A_0=0, 0\leq A_k \leq A_{k+1}$, and let $A(\tau)=\#\{k:T_k \leq \tau\}, \tau\geq 1$. 
 Note that our definition allows more than one event to occur at any given  time point, i.e., batch arrivals. 
%The point arrival process need not be simple (e.g. Bernoulli). 
We also assume that
$A_k \rightarrow \infty$, as $k \rightarrow \infty$, so that there are
only a finite number of arrivals in any finite time interval.  Let $A(\tau):= 
\#\{k: A_k \leq \tau\}$, $D(\tau):= \#\{k: D_k \leq \tau\}, \tau \geq 0$, so that $A(\tau)$ 
and $D(\tau)$ count the number of arrivals and departures, respectively, in the
interval $(0,\tau]$.  Note that, since $A_k < \infty$ for all $k \geq 1$ , $A(\tau)
\rightarrow \infty$ as $\tau \rightarrow  \infty$.  Note also that $A(\tau) = \max
\{k:A_k \leq \tau\}$, since $\{A_k, k \geq 1\}$ is a non-decreasing sequence.
But in general, we cannot write $D(\tau) = \max \{k:D_k \leq \tau\}$, because
$\{D_k, k \geq 1\}$ is not necessarily non-decreasing.  It is
non-decreasing
if the discipline is {\em first-in, first-out} (FIFO), that is, if departures
occur in the same order as arrivals.

Informally, one can see that, 
\bear 
L(\tau) &:=& \#\{k: A_k < \tau \leq D_k\} = A(\tau) - D(\tau) \; , \; \tau \geq 0 \; ,
                                                       \label{eq:defL}      \\
W_k  &:=& D_k - A_k \; , \; k \geq 1 \; ,              \label{eq:defW}
\eear
so that $L(\tau)$ is the number of customers in the system at time $\tau$ and
$W_k$ is the waiting time in the system of customer $k$. 

Formally, let ${\bf 1}\{E\}$ denote the indicator of the event $E$.  Then 
\bear
L(\tau) &=& \sum_{k=1}^\infty {\bf 1}\{A_k < \tau \leq D_k\} \; ,  \label{eq:L}               \\
W_k  &=& \sum_{\tau=0}^\infty {\bf 1}\{A_k < \tau \leq D_k\}   \; .\label{eq:W} 
\eear

 We have seen in Sections~\ref{sec:prel} and \ref{sec:class} that when SR are invoked (\ref{eq:defW}) and (\ref{eq:W}) do not necessarily give the same waiting times $W_k$. It all depends on whether  arrival and departure instants  are based on  observed or actual values. There is no counterpart for this in the continuous time systems.
In this subsection, all we say is that the same arrival and departure instants are used to compute both $W_k$ and $L(\tau)$. 
Our result below is a discrete-time counterpart of the original sample-path version of $L = \lambda W$
contained in Stidham~\cite{Sti72b}, \cite{Sti74}. Other versions, slightly more general, of {\em LL} are given in El-Taha and Stidham~\cite{Elt99}.

\begin{thm} \label{thm:lastword}
Suppose $\tau^{-1}A(\tau) \rightarrow \lambda$ as $\tau \rightarrow  \infty$, where
$0 \leq \lambda < \infty$, and $n^{-1}\sum_{k=1}^n W_k \rightarrow W$
as $n \rightarrow \infty$, where $0 \leq W < \infty$.  Then 
$\tau^{-1}  \sum_{j=1}^\tau L(j) \rightarrow L$, as $\tau \rightarrow  \infty$ and $L = \lambda W$.
\end{thm}

The proof is given in the Appendix. Theorem~\ref{thm:lastword} is general enough to work in almost all practical situations in discrete-time queues. One can  construct  discrete-time versions of {\em LL} that works under weaker conditions than those given here as in El-Taha and Stidham~\cite{Elt99}. Little's law works in discrete-time  with the same interpretation as in the continuous time case. 
%\bigskip

Note that our definitions of $W_k$ and $L(\tau)$ in (\ref{eq:L}) and (\ref{eq:W}) assume $\leq$ on the right  and strict inequality on the left of the indicator function. One can reverse the  equality and inequality without affecting the individual $W_k$ values. However, the $L(\tau)$ values will be affected but the mean value  $L$  will not  as the following example shows. We note that an issue like this does not arise in the continuous time version.
\bigskip

\noindent{\bf Example.} Consider  a discrete system where arrivals and departures occur at  $A_1=1, D_1=4, A_2=2, D_2=5, A_3=5, D_3=7$. The system repeats at time 8, but our focus is in the first  busy period. With the current definition of $W_k$ and $L(\tau)$ in (\ref{eq:L}) and (\ref{eq:W}), we  see that $W_1=3, W_2=3, W_3=2$. 
One can also see that
 $L(1)=0, L(2)=1, L(3)=2, L(4)=2, L(5)=1, L(6)=1$ and $L(7)=1$, so that for the first busy period 
 $L=\tau^{-1}  \sum_{i=1}^\tau L(i)=7^{-1} \sum_{i=1}^7 L(i)=8/7$.

Now, consider  equality on the left and inequality on the right, that is let $W_k = \sum_{\tau=0}^\infty {\bf 1}\{A_k \leq \tau < D_k\}$ and   $L(\tau) = \sum_{k=1}^\infty {\bf 1}\{A_k \leq \tau < D_k\}$. One can see no change in the values of $W_1$, $W_2$ and $W_3$. However,  
$L(1)=1, L(2)=2, L(3)=2, L(4)=1, L(5)=1, L(6)=1$ and $L(7)=0$. Now, for the first busy period 
 $L=\tau^{-1}  \sum_{i=1}^\tau L(i)=7^{-1}  \sum_{i=1}^7 L(i)=8/7$.
In both cases $W=8/3$,  $\lambda=3/7$, $L=8/7$, and  $L=\lambda W$.\done

We see that {\em LL} holds at great level of generality with little assumptions, yet there are concerns that there are discrete-time queueing instances where  {\em LL} appears not to apply. The next  section is devoted to exploring this  apparent contradiction and proposing a resolution. 

\subsection{Applying  Little's Law to  Queues  with Scheduling Rules}

Now let $u(\tau)$ be any observation epoch (see Table~\ref{table-epoch}) and let $A'_k, D'_k, k=1,2,\ldots$ be the {\em SR} dependent arrival and departure times of the $k^{th}$ customer.  Define
\bear
L^o(u(\tau)) &=& \sum_{k=1}^\infty {\bf 1}\{A'_k<u(\tau)\leq  D'_k\} \; ,  \label{eq:Lsr}               \\
W^o_k  &=& \sum_{\tau=0}^\infty {\bf 1}\{A'_k<u(\tau)\leq  D'_k\}   \; ,\label{eq:Wsr} 
\eear
so that $L^o(u(\tau))$ is the observed number of customers in the system at time $u(\tau)$ and
$W^o_k$ is the observed waiting time in the system of customer $k$.
Define the following limits when they exist, 
\bear 
L^o &:=& \limtau\sum_{j=0}^\tau L^o(u(j))/\tau\;,
                                                       \label{eq:defLsr}      \\
W^o  &:=& \limn \sum_{k=1}^n W^o_k/n \; .              \label{eq:defWsr}
\eear

\noindent
{\bf Remarks.}
 When it exists, $L^o$ defined by (\ref{eq:defLsr}) represents the observed long-run time-average number of customers in the system. This time average can be defined  at any of six possible  epochs that we  identified earlier. There may be other epochs of interest. In contrast, in continuous time queues we average continuously over time so there can only be one time-average measure of interest. 

The problem with {\em incoherent} systems is that certain {\em SR}/observation epochs combinations  lead to instances  in which the presence of a customer in the system is not captured by the system state. Capturing  all present customers by the system state is important to identify the correct distribution function and its mean $L$. We explore this further is Section~\ref{sec:LLbdq}. Now, we present our  application of $LL$ to discrete-time queues with scheduling rules.

\begin{thm}\label{thm:LL}
Consider any discrete-time queue where {\em SR} are invoked. Suppose $\tau^{-1}A(\tau) \rightarrow \lambda$ as $\tau \rightarrow  \infty$, where
$0 \leq \lambda < \infty$, and $n^{-1}\sum_{k=1}^n W^o_k \rightarrow W^o$
as $n \rightarrow \infty$, where $0 \leq W^o < \infty$.  Then 
$\tau^{-1}  \sum_{j=0}^\tau L^o(j) \rightarrow L^o$, as $\tau \rightarrow  \infty$ and 
\begin{equation}\label{eq:Lo}
L^o=\lambda W^o\;.
\end{equation}
Moreover, let $W$ be the actual system mean waiting time. Then

(i) for {\em coherent} systems 
\begin{equation}\label{eq:Lcoh}
 L^o=\lambda W\;;
\end{equation}

(ii) for {\em sub-coherent} systems 
\begin{equation}\label{eq:Lsub}
 L^o=\lambda (W-1)\;;
\end{equation}

(iii) for {\em super-coherent} systems 
\begin{equation}\label{eq:Lsup} 
L^o=\lambda (W+1)\;.
\end{equation}
\end{thm}
{\bf Proof.}
Similar to  (\ref{eq:bbasic}), we obtain the inequality,
\beq   \label{eq:bbasic-sr}
\sum_{k:A'_k \leq \tau} W^o_k \geq \sum_{j=1}^\tau L(u(j)) 
                        \geq \sum_{k:A'_k +W^o_k \leq \tau} W^o_k \; , \; \tau \geq 0 \; .
\eeq
The proof of (\ref{eq:bbasic-sr}) is similar to that of Lemma~\ref{lem:ineq}. Now it follows from the definitions, that $W^o_k/k \goto 0$ as $k\goto \infty$, and $k/A'_k \goto \lambda$ as $k\goto \infty$. 
The rest of the proof of (\ref{eq:Lo}) is the same as that of Theorem~\ref{thm:lastword}\;. 
The proof of (i) follows from the definition of $W^o$ for {\em coherent} systems. The proofs of (ii) and (iii) follows from Lemma~\ref{pro:wo}.
\done
\bigskip

 Theorem~\ref{thm:LL} is given at a high level of generality with minimal conditions using sample path analysis. The basic idea of Theorem~\ref{thm:LL} is that for {\em LL} to hold, the measures for the number of customers in the system ($L$) and the time spent in the system ($W$) must be consistently defined.
If you calculate $W$ using the actual time a customer spends in the system, you must also calculate $L$ using the actual number of customers in the system over time.
Similarly, if you calculate $W$ using the observed time (which may differ from the actual time in an incoherent system), you must also calculate $L$ using the observed number of customers at the same observation epochs.

This approach reconciles the apparent contradictions found in the literature. The important thing is to use the same set of arrival and departure instants consistently to evaluate $W$ and $L$. The classification of systems into {\em coherent} and incoherent categories provides a framework for understanding why inconsistent measurements lead to seemingly incorrect applications of the law.
 Little's Law holds universally for all discrete-time queues, regardless of the scheduling rule, as long as this consistency in measurement is maintained. 
\bigskip

\noindent
{\bf Remark.} In (\ref{eq:Lo}), this specialized application  $ L^o=\lambda W^o$ relates the observed $L^o$ to the observed $W^o$ for any discrete-time systems with {\em SR}. It says that if we use the observed values in defining  the waiting times, then we should use the same observed values in computing  the mean queue-length. In  (\ref{eq:Lcoh}), (\ref{eq:Lsub}), and (\ref{eq:Lsup}) we relate the observed $L^o$ to the waiting time, $W$, in the  corresponding actual system. The actual $L$ for the corresponding actual and {\em coherent} systems is given by $\lambda W$. Moreover, $L^o=L$ for {\em coherent} systems, $L^o=L-\lambda$ for sub-coherent systems, and  $L^o=L+\lambda$ for super-coherent systems. 
\bigskip

\noindent
{\bf Remark.} Dattatreya and Singh~\cite{Dat05} study an {\em EAS} model at the random observer and outside observer  epochs. They provide relationships between L values, (their equation (15)), at both epochs using an informal argument. They also  invoke {\em LL} to obtain W, (their  equation (17)). Our approach in Theorem~\ref{thm:LL} provides a formal and rigorous argument  for both equations. Moreover, our results are general enough to include all discrete-time systems with any {\em SR}/observation epochs combinations. 

In the  Section~\ref{sec:LLbdq} below, we show that for queues that can be modeled by birth-death equations, we have one distribution function for each subclass, namely, {\em coherent}, {\em sub-coherent} and {\em super-coherent}  systems. Moreover, we identify the distribution function for each class.

\section{Birth-Death Queues}\label{sec:LLbdq}
 
Here we provide results for discrete-time systems that can be modeled by the birth death equations.
We show that  each class shares the same stationary distribution  regardless of the {\em SR} and observation epochs within its class.

\subsection{General Birth-Death Systems}

We start by considering the birth death equations and give a  general distribution function that we later specialize  for the three classes of queues. 
Recall  the  generalized birth-death (see El-Taha~\cite{Elt23}) equations take the form,
\begin{equation}\label{eq:gbdeqs}
\alpha (n)(1-\beta(n))\pi(n)=(1-\alpha(n+1))\beta(n+1)\pi(n+1), n\geq 0
\end{equation}
 where $\alpha(n)$, $\beta(n)$, and $\pi(n)$  are, respectively,  the state $n$  arrival, service completion, and stationary probabilities.  
Note that the generalized birth-death equations (\ref{eq:gbdeqs}) are valid for any discrete-time queueing system with one potential arrival and one potential departure per time slot.  

To start, we consider a state-dependent Bernoulli queue with state dependent arrival and service completion probabilities.
 Let $\gamma(j)=\frac{\alpha(j)(1-\beta(j)))}{\beta(j+1)(1-\alpha(j+1))}$, $j\geq 0$, and use 
   (\ref{eq:gbdeqs}) to get
\begin{equation}\label{eq:edist2-re}
 \pi(n)= \Pi_{j=0}^{n-1} \gamma(j)\pi(0), n \geq 1
\end{equation}
 where,
\begin{equation}\label{eq:edist3-re}
 \pi(0)= \left[1+ \sum_{k=1}^{\infty} \Pi_{j=0}^{k-1}\gamma(j)\right]^{-1}\;.
\end{equation}

The stationary distribution in (\ref{eq:edist2-re}) is given by Theorem 2.3 of  Daduna~\cite{Dad01} and  El-Taha~\cite{Elt23}. It is, however, important to note that the birth-death equations (\ref{eq:gbdeqs}) and therefore (\ref{eq:edist2-re}) and (\ref{eq:edist3-re}) are valid at any of the six observations epochs.
We use the stationary distribution function in (\ref{eq:edist2-re}) to give  closed form expressions for the  stationary distribution function  for the {\em coherent}, {\em sub-coherent}, and {\em super-coherent} systems.

Consider  the  B/Geom/1 discrete-time queueing system. In this system arrivals follow a Bernoulli process so that  
 $\alpha(n)=\alpha$  for $n\geq 0$. Service time  are geometric with parameter $\beta$. Note that for service completion probabilities $\beta(0)$ and $\beta(1)$ depend on the {\em SR} and observation epochs combination.  Extending the arguments in El-Taha~\cite{Elt23} one can  arrive at the following  service completion probabilities. 
For {\em coherent} systems, we have   $\beta(n)=\beta, n\geq 1$  and $\beta(0)=0$; for {\em sub-coherent} systems, we have  $\beta(n)=\beta,
  n\geq 0$; and  for {\em super-coherent} systems we have  $\beta(n)=\beta, n\geq 2$, $\beta(0)=0$ and $\beta(1)=\frac{\beta(1-\alpha)}{1+\beta}$. 

With this notation  we write $\gamma(0)=\frac{\alpha(1-\beta(0))}{\beta(1)(1-\alpha)}$, $\gamma(1)=\frac{\alpha(1-\beta(1))}{\beta(1-\alpha)}$, and $\gamma(n)=\gamma=\frac{\alpha(1-\beta)}{\beta(1-\alpha)}, n\geq 2$. Simplify to obtain the following result.

\begin{thm}\label{thm:all} Consider  a B/Geom/1 queue as described above. Then
\begin{equation}\label{eq:bd-all}
\pi(n)=\left\{\begin{array}{lll}
        \gamma(0)\gamma(1)\gamma^{n-2}\pi(0)\;;& n=2,3,\ldots\; ;\\\\
        \gamma(0)\pi(0)\;;& n=1\;;
\end{array}  \right.              
\end{equation}
and 
\[
\pi(0)= \frac{1-\gamma}{1-\gamma+\gamma(0)-\gamma(0)\gamma+\gamma(0)\gamma(1)}\;,\]
 where
 for {\em coherent} systems $\gamma(0)=\frac{\alpha}{\beta(1-\alpha)}$, and $\gamma(1)=\gamma$;
 for {\em sub-coherent} systems $\gamma(0)=\frac{\alpha(1-\beta)}{\beta(1-\alpha)}$, and $\gamma(1)=\gamma$;
 and for {\em super-coherent} systems $\gamma(0)=\frac{\alpha(1+\beta)}{\beta(1-\alpha)^2}$, and $\gamma(1)=\frac{\alpha(1+\alpha\beta)}{\beta(1-\alpha)(1+\beta)}$.
\end{thm}

Using Theorem~\ref{thm:all} and simplifying, we compute $L$ for all {\em coherent}, {\em sub-coherent} and {\em super-coherent} systems. 
\begin{eqnarray*}
L= \sum_{n=0}^{\infty}n\pi(n)=\frac{[\gamma(0)\gamma(1-\gamma)^2+\gamma(0)\gamma(1)(1-\gamma(1-\gamma)^2)](1-\gamma)}
{\gamma^(1-\gamma)^2[1-\gamma+\gamma(0)-\gamma(0)\gamma+\gamma(0)\gamma(1)]}\;.
\end{eqnarray*}
We specialize Theorem~\ref{thm:all} for the three classes  in the following result.
\begin{thm} Let $\gamma=\frac{\alpha(1-\beta)}{\beta(1-\alpha)}$ and $\rho=\alpha/\beta$. Then

(i) for {\em coherent} systems
\begin{equation}\label{eq:bd-coh}
\pi(n)=\left\{\begin{array}{lll}
        \rho(1-\gamma)\gamma^{n-1}\;;& n=1, 2,\ldots \;;\\\\
        1-\rho \;;& n=0\;;
\end{array}  \right.              
\end{equation}
(ii) for {\em sub-coherent} systems
\begin{equation}\label{eq:bd-sub}
\pi(n)=\left.  \begin{array}{lll}
        (1-\gamma)\gamma^{n}\;;& n=0,1,\ldots\;;        
\end{array}  \right.              
\end{equation}
(iii) for {\em super-coherent} systems
\begin{equation}\label{eq:bd-super}
\pi(n)=\left\{\begin{array}{lll}
        \rho^2(1-\gamma)\gamma^{n-2}\;;& n=2,3,\ldots \;;\\\\
        (\alpha+\rho)(1-\rho)\;;& n=1 \;;\\\\
         (1-\alpha)(1-\rho)\;;& n=0\;.
\end{array}  \right.              
\end{equation}
\end{thm}

This theorem generalizes known results to all subsystems. It shows that all {\em coherent} birth death queues share the   same distribution function. Similarly, all {\em sub-coherent}, and {\em super-coherent} systems share  similar distribution functions.  {\em Coherent} systems are the ones of primary interest. Other systems can be of interest in certain cases. El-Taha~\cite{Elt24} shows that  in the birth death  {\em EAS} and {\em LA-DF} models with Bernoulli arrivals,  the distribution function at pre-arrival  epochs  coincides with the distribution function of the {\em sub-coherent} birth death  systems.
\bigskip

\noindent
\subsection{Coherent and Sub-Coherent Systems}

It is more common in the literature to see  studies of instances of {\em coherent} and/or {\em sub-coherent} systems as the majority of systems
 fall under these two categories. For these systems $\gamma(1)=\gamma$, and $\gamma(0)=\frac{\alpha(1-\beta(0))}{\beta(1-\alpha)}$ where $\beta(0)=0$ for {\em coherent} systems and $\beta$ for {\em sub-coherent} systems. Therefore,
\begin{equation}\label{eq:bd-gbd1}
\pi(n)=\left\{\begin{array}{lll}
        \frac{\gamma(0)(1-\gamma)\gamma^{n-1}}{1-\gamma+\gamma(0)}\;;& n=1,2,\ldots \;;\\\\
        \frac{1-\gamma}{1-\gamma+\gamma(0)}\;;& n=0\;;
\end{array}  \right.              
\end{equation}
 and
\[
L= \sum_{n=0}^{\infty}n\pi(n)=\frac{(1-\gamma)\gamma(0)}{1-\gamma+\gamma(0)}\sum_{n=0}^{\infty}n\gamma^{n-1}
=\frac{(1-\gamma)\gamma(0)}{1-\gamma+\gamma(0)}\frac{1}{(1-\gamma)^2}\;;
\]
then simplify to obtain,
\[
L =\frac{\gamma(0)}{(1-\gamma+\gamma(0))(1-\gamma)}\;.
\]
Using the input parameters $\alpha$ and $\beta$, we get
\[
L=\frac{\alpha-\alpha\beta(0)}{\beta-\alpha\beta(0)}\times\frac{\beta(1-\alpha)}{\beta-\alpha}\;;
\]
so that 
\begin{equation}\label{eq:bd-L}
L=\left\{\begin{array}{lll}
    \frac{\alpha(1-\alpha)}{\beta-\alpha}\;;& \beta(0)=0\; ($coherent systems$)\;; \\\\
     \frac{\alpha(1-\beta)}{\beta-\alpha}\;;& \beta(0)=\beta\; ($sub-coherent systems$)\;.
\end{array}  \right.              
\end{equation}

Now we focus on $\pi(0)$. Note that $1-\gamma+\gamma(0)=\frac{\beta-\alpha\beta(0)}{\beta(1-\alpha)}$, and 

\[ \pi(0)=\frac{\beta-\alpha}{\beta-\alpha\beta(0)}\;. \]

Moreover, the probability that there is one or more customers in the system is given by
\begin{equation}\label{eq:util}
1-\pi(0)=\frac{\gamma(0)}{1-\gamma+\gamma(0)}=\frac{\alpha-\alpha\beta(0)}{\beta-\alpha\beta(0)}\;.
\end{equation}
Equation~(\ref{eq:util}) is interesting in that  one cannot say, at this level of generality, that $1-\pi(0)$ is the probability the server is busy, or it represents the sever utilization, only that it represents the probability that one or more customers are in the system. This is a radical shift from our understanding of
$1-\pi(0)$ is continuous time systems.  
Note how we adjust for customers that enter an empty system and leave immediately, i.e., before the next observation instant.
\bigskip

\noindent
{\bf Remark.}
In the continuous setting we apply Little's law to obtain the server's utilization where $L=U$, $\lambda$ is the arrival rate and $W$ is the mean service time. Applying {\em LL} to the discrete-time {\em actual} system, we obtain
$
U=\alpha/\beta
$
as the probability that the server is busy.  This is true only for {\em coherent} systems.
Similarly, the probability the server is idle is $I=1-\rho$. Now, note  that,
\begin{equation}\label{eq:bd-gbd2b0}
\pi(0)=\left\{\begin{array}{lll}
	\frac{\beta-\alpha}{\beta}=1-\rho\;;& \beta(0)=0\; ($coherent systems$)\;; \\\\
         \frac{\beta-\alpha}{\beta-\alpha\beta(0)}=1-\gamma\;;& \beta(0)=\beta\; ($sub-coherent systems$) \;.
\end{array}  \right.              
\end{equation}
So, the time average probability that the server is idle is consistent with the case $\beta(0)=0$. For {\em sub-coherent} systems where
 $\beta(0)=\beta$,  we cannot interpret $\pi(0)$  as the probability the server is idle, only that there are $0$ customers in the system at the observed epochs.  This is an adjustment  for customers that enter an empty system and leave immediately, i.e., before the next observation instant. Obviously, those customers kept the server busy for one slot, but that information is not captured in the state description.

%%%%%%%%%%%%
\section{Additional Invariant Characteristics of Coherent Systems}\label{sec:invar}
%%%%%%%%%%%%%

In this Section we focus on {\em coherent} systems and give results that are invariant with respect to {\em SR}.  Note that {\em incoherent} systems lead to inconsistencies in the sense that the time a server spends busy with a customer is not exactly the same as the service requirement of the said customer. {\em Coherent} systems are more internally consistent and more widely used in the literature.   
In the next two subsections we focus on busy servers and busy periods. 

\subsection{Busy Servers}

Consider a stable  multi-server discrete-time queueing system with $c$ parallel identical servers, and general inter-arrival times and service times distribution functions, and a non-preemptive queueing discipline.  Since no tracking of the  system's state is necessary,  we  use the {\em actual} system   without specifying any scheduling rules. Our first result concerns the mean number of busy servers in a stable multi server system.

Let $\{A_n,S_n, n\geq 1\}$ be the input data for this system where $A_n$ is the arrival instant of the $n^{th}$ arrival and $S_n$ is the service requirement of the $n^{th}$ arrival. Assume that for all $n$, $A_n\leq A_{n+1}$, that is we allow batch arrivals. Also let $A_n \goto \infty$ as $n\goto\infty$.
Let $A(\tau)=\max\{n:A_n\leq \tau\}$ be the number of arrivals during $(0,\tau]$.
  Define the following limits when they exist,
\[
\alpha=\limtau A(\tau)/\tau=\limn n/A_n\;;
\]
\[
ES=\limn \sum_{k=1}^n S_k/n\;.
\]

For $i=1,\ldots,c$, let  $B_i(\tau)=1$ if the $i^{th}$ server is busy at $\tau$, and $0$ otherwise. Also let $Y_i(B,\tau)=\sum_{u=1}^{\tau}\bo\{B_i(u)=1\}$
 be the total time the $i^{th}$ server is busy during $(0,\tau]$.
Define the long-run fraction of time the $i^{th}$ server is busy as
\[
U_i =\limtau Y_i(B;\tau)/\tau\;.
\]
Then  $U=\sum_{i=1}^{c}U_i$,  the sum of the servers' busy fractions, can be interpreted as the long run average number of busy servers. 

\begin{lem}\label{lem:U} Consider a stable ($\rho:=\frac{\alpha ES}{c}<1$) non-preemptive work conserving $GI/GI/c$ discrete-time queueing system. At the start of service, a customer that finds  idle servers, will pick one at random.
 Let 
 $U=\sum_{i=1}^{c}U_i$, then for {\em coherent} systems
\[
U=\alpha ES
\]
\end{lem}
{\bf Proof.} For $i=1,\ldots,c$, let  $S^i_n$ be the service requirement of $n^{th}$ departing customer, and $D^i_n$ is the departure time of the $n^{th}$ departing customer.
Then  the busy time of the $i^{th}$ server during $(0,D^i_n]$  is equal to sum of service completions, i.e.,  
$
Y_i(B,D^i_n)=\sum_{k=1}^nS^i_k\;.
$ 
Let $\alpha^i$ and $\delta^i$ be the rate  at which customers, respectively, join and depart server $i$.  Now

\begin{eqnarray*}
U&=&\sum_{i=1}^cU_i\\
&=&\sum_{i=1}^c \limtau\frac{Y_i(B, \tau)}{\tau}\\
%&=& \limtau \sum_{i=1}^c\frac{Y(i, \tau)}{\tau}\\
&=& \sum_{i=1}^c \limn \frac{Y_i(B,D^i_n)}{D^i_n}\\
&=& \sum_{i=1}^c \limn \frac{\sum_{k=1}^{n}S^i_k}{D^i_n}\\
&=& \sum_{i=1}^c \limn \frac{n}{D^i_n}\frac{\sum_{k=1}^{n}S^i_k}{n}\\
&=& \sum_{i=1}^c \alpha^iES^i\\
&=& \alpha ES\,.
\end{eqnarray*}
Here we used  the coherence assumption, a discrete version of Lemma 2.1 of El-Taha~\cite{Elt99}, $ES^i=ES$, $\delta^i=\alpha^i$ and $\alpha=\sum_{i=1}^c\alpha^i$. This completes the proof.\done

Because of Lemma~\ref{lem:U}, one can assert that $L=L_q+\lambda ES$, where $L_q$ is the mean number of customers in the queue (excluding time in service). This known relationship applies for {\em coherent} systems only. Now, we focus on single server systems.

\begin{cor} Consider a stable ($\rho:=\alpha ES<1$) $GI/GI/1$ discrete-time {\em coherent} queueing system. Let 
\[
U =\limtau Y(B;\tau)/\tau
\]
where $Y(B,\tau)$ be the total time the server is busy during $(0,\tau]$. Assume that $U$ is well-defined, then,
\[
U=\alpha ES\;.
\]
\end{cor}
{\bf Proof.} The proof follows from Lemma~\ref{lem:U}.\done

\begin{cor}
For {\em coherent}  stable single server GI/GI/1 systems, let $\pi(0)$ be the  stationary probability that the server is idle, i.e., there are $0$ customers in the system. Then
\[
1-\pi(0)=\rho \;.
\]
\end{cor}
{\bf Proof.} $1-\pi(0)=\limtau \tau^{-1}\sum_{k=1}^{\tau}\bo\{Z(k)>0\}$ w.p.1. That is the long-run fraction of time the server is busy. By Lemma~\ref{lem:U} this is equal to $\alpha ES$.\done

\noindent
{\bf Remark.} One can obtain $U$ in Lemma~\ref{lem:U} for {\em coherent} systems using {\em LL}, where 
the system (black box) is the service system itself, so that $L$ is the mean number of busy servers and $W$ is the mean service time.

Now, let $I$ be the long-run fraction of  time that the server is idle, then $I=1-U$, so that
\[
I=1-\alpha ES \;.
\]
In single server queues, $I$ represent the long-run fraction of time the server is idle. In continuous time queues this is the same as $p(0)$, the probability that there are  $0$ customers in the system. For discrete-time queues this is true only for {\em coherent} systems.
 That is, it  depends on the {\em SR} and the observation instants as shown in the following example.
 \bigskip

\noindent
{\bf Example.} Consider the B/Geo/1 queue and compute $1-\pi(0)$ the probability that there is at least one customer in the system, for all combinations of {\em SR}/observation epochs. The results are given in Table~\ref{table-exp} below. 

\begin{center}
\noindent
\begin{tabular}{|c|ccccccc|}
\multicolumn{6}{c}{\bf Table~\ref{table-exp}. $1-\pi(0)$  for single server Markovian queues} \\
\multicolumn{4}{c}{} \\
 \hline  
 & Random & Outside & Poten.-Pre& Poten.-Post&Poten.-Pre&Poten.-Post  &\\
& Observer& Observer& Arrival&Arrival&Departure&Departure&\\ 
\hline 
&&&&&&&\\  
$EAS$       &$\gamma$ & $\rho$ &$\gamma$&$\rho$ &$\rho$ &$\gamma$& \\&&&&&&& \\
LAS-IA      &$\rho$ &$\gamma$&$\gamma$ &$\rho$ &$\rho$ &$\gamma$&\\ &&&&&&&\\
LAS-DA   &$\rho+$ &$\rho$  & $\rho$ &$\rho+$&$\rho+$&$\rho$ &\\&$\alpha(1-\rho)$&&&$\alpha(1-\rho)$&$\alpha(1-\rho)$&& \\
&&&&&&&\\
LA-AF   &$\rho$ &$\rho$  & $\rho$ &$\rho+$&$\rho+$&$\rho$ &\\&&&&$\alpha(1-\rho)$&$\alpha(1-\rho)$&& \\
&&&&&&&\\
 LA-DF  &$\rho$ &$\rho$ &$\gamma$&$\rho$ &$\rho$ &$\gamma$&\\ &&&&&&&\\
\hline
\end{tabular}
\label{table-exp}
\end{center}
\vspace{5mm}

Note that in seventeen out thirty possible combinations, we obtain  $1-\pi(0)=\rho$ as expected. These combinations represent the {\em coherent} discrete-time queues. In all other cases we obtain a value for $1-\pi(0)=\gamma$ for {\em sub-coherent} systems and $1-\pi(0)=\rho+\alpha(1-\rho)$ for {\em super-coherent} systems. For {\em incoherent} systems $1-\pi(0)$ does not represent the true value of the server's utilization factor. Using {\em LL} for sub-coherent systems, we see that $U=\alpha (ES-1)=\rho(1-\beta)=\gamma(1-\alpha)$. On the other hand, using {\em LL} for super-coherent systems,  $U=\alpha(ES+1)=\rho(1+\beta)$.

\subsection{Busy Periods}

In this section we provide busy period analysis using the actual  discrete-time systems as described in Section~\ref{sec:LL}, so that the results will be valid for all {\em coherent} systems.
We assume work conserving queueing discipline, i.e., the server is not idle when there is work in the system.  
Let  the random variable $A$ represent inter-arrival times and the random variable $S$ represent service times.
 Assume that inter-arrival times and service times are $i.i.d.$ and  independent of each other. Let the mean inter-arrival times 
$E(A)=1/\alpha$, and mean service times $ES=1/\beta$  where $0<\alpha,\beta<\infty$, and let  the traffic intensity $\rho=\alpha/\beta <1$. Assume the system is empty at time $0$.
For $j=1,\ldots$, let
\[
a(u(j))=\bo \{L(u(j)-1)=0,L(u(j))\geq 1\}\;,
\]
and
\[
d(u(j))=\bo \{L(u(j)-1)\geq 1,L(u(j))=0\}\;.
\] 
Then $A_{B1}=\{u(j): a(u(j))=1\}$ represents the set of arrival instants that find the system idle for at least one unit, i.e., the set of busy period start instants. Moreover, $D_{B1}=\{u(j):d(u(j))=1\}$ is the set of departure instants that leave the system idle after being busy for at least one period, i.e., idle period start instants.

Let 
$A_B(0,\tau)=\sum_{j=1}^{\tau} a(u(j))$  and $D_B(0,\tau)=\sum_{j=1}^{\tau} d(u(j))$  be,  respectively,  the number of arrival and departure instants to find (leave) the system empty during $(0,\tau ]$.
It follows from the definitions that for all $\tau \geq0$,  $0\leq A_B(0,\tau)-D_B(0,\tau)\leq 1$.
Moreover, let
$U_k=\min\{\tau: A_B(0,\tau)=k\}$ and $V_k=\min\{\tau: D_B(0,\tau)=k\}$ be the $k^{th}$ arrival (departure) to find (leave) the system idle. Note that for all $k=1,2,\ldots,$, $U_k<V_k<U_{k+1}$. Let
\begin{eqnarray*}
C_k&=&U_{k+1}-U_k\;;\\
B_k&=&V_k-U_k\;;\\
I_k&=&U_{k+1}-V_{k}\;;\\
E_k&=&A(U_{k+1})-A(U_{k})\;.
\end{eqnarray*}

We interpret $C_k$ as the length of the $k^{th}$ busy cycle, $B_k$ the length of the $k^{th}$ busy period, $I_k$ as the length of the $k^{th}$ idle period and $E_k$ as the total number of arrivals (service completions) in the $k^{th}$ busy period (cycle). Note that $C_k=B_k+I_k$.

 Define the following limits when they exist,
\begin{eqnarray}
I&=&\limn \sum_{k=1}^n I_k/n\;;\\
C&=&\limn \sum_{k=1}^n C_k/n\;;\\
B&=&\limn \sum_{k=1}^n B_k/n\;;\\
E&=&\limn \sum_{k=1}^n E_k/n\;.
\end{eqnarray}

Our results  are valid for the actual/standard queueing model. 
Note that for {\em coherent} systems, the observed and actual service times are the same, and that  the busy period is the sum of all service times served during a busy period, so the busy period is the same for all {\em coherent} systems. Server's busy times are invariant  for {\em coherent} systems. This shows that our results in this section are valid for all {\em coherent} systems.
The following result follows immediately from these observations.
\begin{lem}\label{lem:busyP}
For all $k=1,\ldots$, $B_k$, $C_k$, $I_k,$ and  $E_k$ are invariant with respect to {\em coherent} systems and equal to the corresponding actual system quantities.
\end{lem}

Let
$Y(n,\tau)=\sum_{j=1}^{\tau}\bo\{L(u(j))=n\}$ be the time spent in state $n$ during $(0,\tau]$, and $A(n,\tau)$  be the number of arrivals  that find  process $Z$ in state $n$  during $(0,\tau]$. Define the following limits when they exist.
\begin{eqnarray}
\pi(n)&=&\limtau Y(n;\tau)/\tau\;;\\
\alpha(n)&=&\limtau A(n;\tau)/ Y(n;\tau)\;.
%\beta(n)&=&\limtau D(n;\tau)/ Y(n;\tau)\;.\\
\end{eqnarray}
That is $\pi(n)$ is long-run fraction of  time process $\{L(\tau), \tau>0\}$ is in state $n$, and $\alpha(n)$ is the state-n arrival rate.
Now, we state the following result.
\begin{thm}\label{thm:busyperiod} 
Assuming the limits exist, then 
\begin{eqnarray}
I&=&1/\alpha(0)\;;\\
C&=&1/\alpha(0)\pi(0)\;;\\
B&=&(1-\pi(0))/\alpha(0)\pi(0)\;;\\
E&=&\alpha/\alpha(0)\pi(0)\;.
\end{eqnarray}
\end{thm}
{\bf Proof.} Using $Y=\lambda X$ (El-Taha and Stidham~\cite{Elt99}), we have
\begin{eqnarray}
I&=&\limtau Y(0;\tau)/A(0;\tau)\;;\\
B&=&\limtau (\tau-Y(0;\tau))/A(0;\tau)\;;\\
C&=&\limtau \tau/A(0,\tau)\;;\\
E&=&\limtau A(\tau)/A(0;\tau)\;.
\end{eqnarray}
We show how $Y=\lambda X$ is used to prove first case.
Let $Y(\tau)=Y(0;\tau)$, and $X_k=Y(0,U_{k+1})-Y(0,U_{k})=I_{k}$, $k\geq 1$. Here $\lambda =\limtau A(0,\tau)/\tau$.
Now,
\begin{eqnarray} 
I&=& X=Y/\lambda=\limtau  (Y(0;\tau)/\tau) (\tau/A(0,\tau) \\
&=& \limtau  Y(0;\tau)/A(0,\tau)\;.
\end{eqnarray}
Therefore $I=1/\alpha(0)$. The proof of the other results is similar.

%%%%%%%%%%%%
\subsubsection{Applications}
%%%%%%%%%%%%%

Here we give three examples where using input parameters we are able to compute the mean  busy and idle periods, mean busy cycles, and the mean number of arrivals during a busy period. 
\bigskip

\noindent
{\bf Example 1.} Consider a stable $B/G/1$ {\em coherent} queueing model, i.e., $\rho=\alpha/\beta <1$. Then it follows  that 
$\pi(0)=1-\rho$. Now, because arrivals are Bernoulli and the arrival process and service times are independent,  using {\em BASTA} (e.g., El-Taha~\cite{Elt24}), we have   $\alpha(0)=\alpha$ a.s.. Thus
\begin{eqnarray*}
I&=&1/\alpha\;;\\
C&=&1/\alpha(1-\rho)\;;\\
B&=&1/\beta(1-\rho)\;;\\
E&=&1/(1-\rho)\;.
\end{eqnarray*}

\bigskip

\noindent
{\bf Example 2.} 
Consider  a {\em coherent}  discrete-time finite population model  $B/Geom/1//N$. Here  $\alpha (0)=N\alpha$, $\alpha_{eff}=\alpha(N-L)$ where $L$ is the mean number of customers in the system, and $\pi(0)$ is computed numerically as in El-Taha~\cite{Elt23}. Thus
\begin{eqnarray*}
I&=&1/N\alpha\;;\\
C&=&1/N\alpha \pi(0)\;;\\
B&=&(1-\pi(0))/N \alpha \pi(0)\;;\\
E&=&(N-L)/N \pi(0)\;.
\end{eqnarray*}

\noindent
{\bf Example 3.} Consider a discrete-time  $G/Geo/1$ {\em coherent}-queue. We need the distribution function at the pre-arrival and potential pre-arrival instants. Note that only the {\em LAS-DA} and {\em LA-AF} {\em SR} are {\em coherent} at  pre-arrival and potential pre-arrival epochs.

Let $\pi(.)$ and $\pi^A(.)$ be, respectively, the potential pre-arrival  and pre-arrival probabilities.
 Here $\pi(0)=1-\rho$. Note that $\pi^A(0)=1-\sigma^* $ where $\sigma^*=\sigma/(\sigma\beta+1-\beta)$ and $\sigma$ is the unique solution in $(0,1)$ of $\sigma=F^*(\sigma\beta+1-\beta)$, $F(.)$ is the inter-arrival time distribution function, and $F^*(.)$ is its probability generating function.  See Hunter~\cite{Hun83b}, pp 449-251 for details.
 We also know that $\alpha(0)\pi(0)=\alpha\pi^A(0)$. 
 Thus $\alpha(0)= \alpha\pi^A(0)/\pi(0)=\alpha(1-\sigma^*)/(1-\rho)$. 
Therefore
\begin{eqnarray*}
I&=&(1-\rho)/\alpha(1-\sigma^*)\;;\\
C&=&1/\alpha (1- \sigma^*)\;;\\
B&=&\rho/\alpha (1- \sigma^*)\;;\\
E&=&1/(1- \sigma^*)\;.\done
\end{eqnarray*}

One can apply Theorem~\ref{thm:busyperiod} to other examples as well.
In the next section we use $H=\lambda G$ to give a proof for the PK formula that is valid for all {\em coherent} systems.

\section{Discrete-Time H=$\lambda$G and Waiting Times}\label{sec:HlG}

Consider the deterministic sequence of time points $\{A_k, k\geq 1\}$ as given in  Section~\ref{sec:LL}.
Associated with each time point $A_k$, there is a function  
$f_k: I \rightarrow R^+$, where $I$ is the set of non-negative integers and $R^+$ is the set of non-negative real numbers. We assume that  $f_k(\tau)$  is Lebesgue integrable on  
$\tau \in [0,\infty)$, for each $k \geq 1$.
The bivariate sequence  $\{(A_k,f_k(\cdot)), k \geq 1\}$  
constitutes the basic data, in terms of which the behavior of the system is 
described.
Let $f_k(\tau)$ denote the rate at which customer $k$ incurs cost at time $\tau$,
$k \geq 1$, $\tau \geq 0$.  Define
\bear  
H(\tau) &:=& \sum_{k=1}^{\infty} f_k(\tau) \; , \; \tau \geq 0 \; , \label{eq:defHt}  \\
G_k  &:=& \sum_{j=0}^{\infty} f_k(j)   \; , \; k \geq 1 \; , \label{eq:defGk}
\eear
so that $H(\tau)$ is the total cost rate at time $\tau$ and $G_k$ is the total cost 
incurred by customer $k$. Assume that $H(\tau)$ and $G_k$ are well defined for all $\tau$ and $k$.

Like the continuous time case, {\em LL} has an economic interpretation that  suggests the current extension. Let $f_k(\tau): = {\bf 1}\{A_k < \tau \leq D_k\}$, i.e., customer $k$ incurs a cost of one dollar per unit time while
in the system (i.e., while $A_k < \tau \leq D_k$) and zero cost otherwise.
 Then we can interpret the 
function $f_k(\tau)$ as the cost rate of customer $k$ at time $\tau$.  Under this 
interpretation, 
$L(\tau) = \sum_{k=1}^{\infty} f_k(\tau)$ is the total cost rate at time $\tau$ and
$W_k = \sum_{\tau=0}^\infty f_k(\tau) $ is the total cost incurred by customer $k$, so that 
$L = \lambda W$ 
says that the long-run average cost per unit time 
equals the arrival rate of customers times the long-run average cost per
customer.
The generalization to $H = \lambda G$ arises
 naturally if one allows a
more general cost-rate function than the indicator of the event 
$\{A_k < \tau \leq D_k\}$.  
With $H(\tau)$ and $G_k$ defined by (\ref{eq:defHt}) and (\ref{eq:defGk}), 
respectively, define the following limiting averages, when they exist:
\begin{eqnarray}
\lambda       &:=&  \limtau \tau^{-1}  N(\tau) \; ,              \label{deflam2}   \\
H             &:=&  \limtau \tau^{-1}  \sum_{j=0}^\tau H(j) \; ,      \label{defH}      \\
G             &:=&  \limn n^{-1} \sum_{k=1}^n G_k  \; . \label{defG}
\end{eqnarray}

We seek conditions under which $H = \lambda G$.  
Following Stidham~\cite{Sti74} and Heyman and Stidham~\cite{Hey80}, suppose 
that the bivariate sequence $\{(A_k,f_k(\cdot)), k \geq 1\}$ satisfies the 
following condition:
\bigskip

\noindent
{\bf Condition L1.}  There exists a sequence $\{W_k, k \geq 1\}$ such that,

(i)  $W_k/A_k \rightarrow 0$ as $k \rightarrow \infty$ ; and

(ii) $f_k(\tau) = 0$ for $\tau \notin (A_k,A_k+W_k]$ .
\bigskip

Condition {\em L1} says that all the cost associated with the 
$k^{th}$ point (e.g., the $k^{th}$ customer) is incurred in a finite time 
interval beginning at the point (e.g., the arrival of the customer), and that 
the lengths of these intervals cannot grow at the same rate as the points 
themselves, as $k \to \infty$.  This is a stronger-than-necessary condition 
for $H = \lambda G$ (See El-Taha and Stidham~\cite{Elt99} for details), but it is 
satisfied in most applications to queueing systems, in which the time points 
$A_k$ and $A_k+W_k$ correspond to customer arrivals and departures, 
respectively, and it is natural to assume that customers can only incur cost 
while they are physically present in the system.   

The proof of the discrete-time $H=\lambda G$ follows the same steps as the continuous-time case given by El-Taha and Stidham~\cite{Elt99}.

\begin{thm} \label{thm:pmtHlG}
Suppose $\tau^{-1} A(\tau) \rightarrow \lambda$ as $\tau \rightarrow  \infty$, where
$0 \leq \lambda < \infty$, and Condition L1 holds.  Then if $n^{-1}\sum_{k=1}^n G_k \rightarrow G$ as $n \rightarrow \infty$, 
where $0 \leq G < \infty$, then $\tau^{-1}  \sum_{j=0}^\tau H(j)$ $\rightarrow H$ as 
$\tau \rightarrow  \infty$, and $H = \lambda G$.
\end{thm}
{\bf Proof.} Similar to (\ref{eq:bbasic}), one can show that for all $\tau \geq 1,$
\beq   \label{eq:hlgbasic}
\sum_{k:A_k \leq \tau} G_k \geq \sum_{j=1}^\tau H(j) 
                        \geq \sum_{k:A_k+W_k \leq \tau} G_k \; , \; \tau \geq 1 \; .
\eeq
Now, utilizing Lemma~\ref{lem:nt}, the proof is similar to the argument in the proof of Theorem~\ref{thm:lastword}.\done

\noindent
\subsection{Multi-Server Queues: Relations between Workload and 
Waiting Time}\label{subsection:ggc}

The results in this Subsection are derived using the {\em actual} system, i.e., we assume  arrivals and departures take place at  the integer instants $\tau$.  By condition $L2$ (see below), the results are valid for all {\em coherent} systems. 
We now  use $H = \lambda G$ to derive a relation between the
time-average workload and the customer-average waiting time in the queue in
a multi-server system with a non-preemptive queue discipline. 
Consider the {\em G/G/c} queue.  The 
input data consists of the sequence 
$\{(A_k,S_k), k \geq 1\}$, where $A_k$ is the arrival instant and $S_k$ the 
work requirement of customer $k$.  
 Assume that each the queue discipline is a non-preemptive, i.e.,
the server is never idle 
when customers are waiting, and that the server works at unit rate.

Let
$W_k^q$ denotes the $k^{th}$ customer waiting time in queue (excluding service time).
Assume  the following limits exist and are finite:
\bears
ES      &:=& \limn n^{-1} \sum_{k=1}^n S_k  \; ,    \\
ES^2    &:=& \limn n^{-1} \sum_{k=1}^n S_k^2  \; ,  \\
EW^q    &:=& \limn n^{-1} \sum_{k=1}^n W_k^q \; , \\
ESW^q   &:=& \limn n^{-1} \sum_{k=1}^n S_k W_k^q \; .
\eears
Here, $ES$ is the long-run average service time, $ES^2$ is the long-run  empirical second moment of service times, $EW_q$ is the long-run average waiting time in queue. Note that these are sample-path averages, even though we use a notation
suggestive of expectations. 

Let  $f_k(\tau)$ be the work remaining to be done for the $k^{th}$
customer at time epoch $\tau$.  Then for the {\em actual} system
\bear \label{eq:fktau}
f_k(\tau) &=& S_k {\bf 1}\{ A_k < \tau \leq A_k+W_k^q \}  +     \nonumber           \\
       & &   (S_k-(\tau-A_k-W_k^q))
               {\bf 1}\{ A_k+W_k^q < \tau \leq A_k+W_k^q+S_k \} \; .
\end{eqnarray}

Implicit in (\ref{eq:fktau}) is the assumption that all arrivals and departures take place at the discrete-time instants $\tau$.
Let
\[
V(\tau)= \sum_{k=1}^{\infty} f_k(\tau),
\]
so that  $V(\tau)$ is the total amount of work in the system at potential arrival times $\tau$. That is this is the total work in the system a potential arrival at time $\tau$ would see upon arrival.
Let
\[
EV:= \limtau \tau^{-1}  \sum_{j=1}^\tau V(j)  \; ,
\]
when the limit exists.  
The following result is valid for  multi-server discrete-time queues.

\begin{cor} Consider a multi-server queue with {\em FIFO} queueing discipline. Suppose the sequences $\{S_k, k \geq 1\}$ and $\{W_k^q, k \geq 1\}$, are
{\em asymptotically pathwise uncorrelated}, that is,

\beq   \label{eq:apunc}
ESW^q = \limn n^{-1} \sum_{k=1}^n S_k W_k^q = ES \cdot EW^q \; .
\eeq
Then
\begin{equation}  \label{eq:EV}
EV = \lambda ES EW^q + \lambda (ES^2-ES)/2 \; .
\end{equation}
\end{cor}
{\bf Proof.}
Let $H(\tau) = V(\tau)$. Using (\ref{eq:fktau}),
\[
G_k = \sum_{\tau=0}^\infty f_k(\tau)= S_k W_k^q + (S_k^2-S_k)/2 \; ;
\]
\[
G = \limn n^{-1}\sum_{k=1}^{n}[S_k W_k^q + (S_k^2-S_k)/2] = ESW^q + (ES^2-ES)/2 \; .
\]
Since $\lambda$, $ES$, and $EW^q$ are well defined
and finite, Condition $L1$ holds with $W_k = W_k^q + S_k$, the waiting time
of the $k^{th}$ customer in the system.  Applying $H = \lambda G$, we conclude that

\begin{equation}  \label{eq:avgwkld}
EV = \lambda ESW^q + \lambda (ES^2-ES)/2 \; .
\end{equation}
Now, use condition (\ref{eq:apunc}) to obtain (\ref{eq:EV}).\done

Condition (\ref{eq:apunc}) is  true a.s. for stochastic models with {\em service-time
independent} queueing disciplines, that is, models in which the rule for selecting
the next job to process does depend on the service times
of jobs.  The {\em FIFO} queue discipline is an example of such a rule.

The first term of (\ref{eq:EV}) is the total amount of work associated with 
customers waiting in the queue, and the second term is the 
residual service time. In contrast, for  continuous time models the residual service time is given by $\lambda ES^2/2$.
Relationship (\ref{eq:EV}) is valid for systems where customers depart in the order of arrival like the   $G/D/c$ and $G/G/1-FIFO$ queues.

Note that $EV$ 
coincides with the
{\em virtual} waiting time, that is, the amount of time that a customer would have
to wait in the queue if that customer arrived at time $\tau$.  Thus, in the
{\em FIFO} case, (\ref{eq:EV}) also gives a relation between the 
time average virtual
waiting time and the customer average actual waiting time in queue.
\bigskip

\noindent
{\bf Single-Server Queues: The Pollaczek-Khinchine Formula}
\medskip

Here we give the well-known Pollaczek-Khinchine Formula for the {\em actual} discrete-time queues with Bernoulli arrivals.

\begin{cor}\label{cor:spPK} Consider a stable {\em FIFO} single-server queue with service time independent  discipline.
Assume $ASTA$ holds, then,
\beq  \label{eq:spPK}
EW^q = \frac{\lambda (ES^2-ES)}{2(1 - \rho)} \; ,
\eeq
where $\rho= \lambda ES <1$.
\end{cor}
{\bf Proof.} 
Let $EV^A$ be  the workload in the system at (actual) arrival times. 
 With {\em FIFO} discipline, the workload at arrival instants $EV^A$ coincides with the actual waiting time in the queue $W^{q}$.
Moreover, by {\em ASTA} (El-Taha~\cite{Elt24}), we obtain $EV^A=EV$ (note that EV is the workload at {\em potential} arrival time instants). 
Therefore
\[
EW^{q} = \lambda ES EW^{q} + \lambda (ES^2-ES)/2 \;.
\]
Simplify to obtain (\ref{eq:spPK}).\done

 Corollary~\ref{cor:spPK} is a 
sample-path version of the Pollaczek-Khinchine (PK) formula for a {\em FIFO} single-server
queue where {\em ASTA} holds.  In a stochastic setting, the most common situation where   {\em ASTA} holds is when arrivals are Bernoulli.

\subsection{H=$\lambda$G and Coherent Queues}

Here we explore how H=$\lambda$G applies to queues with {\em SR}. We first extend our definition of {\em coherent} systems to  queues with cost function $f_k(\tau)$.
Let $A_k,D_k=A_k+W_k$ be the actual system $k^{th}$ arrival and departure epochs. Moreover, let $A'_k, D'_k=(A_k+W_k)'$ be the $k^{th}$  scheduled arrival and departure epochs associated with selected  {\em SR}. Now, we modify condition {\em L1} to {\em L2}. 

\noindent
{\bf Condition L2.} 
There exists a sequence $\{W_k, k \geq 1\}$ such that,

(i)  $W_k/A'_k \rightarrow 0$ as $k \rightarrow \infty$ ; and
%\[

(ii) $f_k(u(\tau))=0 \mbox{ for } u(\tau) \notin (A'_k, (A_k+W_k)']$
%\]
where depending on the  {\em SR} we always have $A'_k$ equals  $A_k^-, A_k^{--}$ or $A_k^+$. Similarly,  $D'_k=(A_k+W_k)'$ equals  $D_k^-, D_k^{--}$, $D_k^+$ or $(D_k-1)^+$.

\noindent
{\bf Definition 1.} Consider  a discrete-time queueing system with any {\em SR} and observation epoch combination.  
Then, the system is said to be {\em coherent} if for all arrivals $k=1,2,\ldots $,
\[
\sum_{\tau=1}^{\infty}f_k(\tau)= \sum_{\tau=1}^{\infty}f_k(u(\tau))\;.
\]
Otherwise, the system is said to be {\em incoherent}. \done

Note that for all $k=1,\ldots $
$
\sum_{\tau=1}^{\infty}f_k(\tau)=\sum_{\tau=A_k+1}^{A_k+W_k}f_k(\tau), %(D_k=A_k+W_k)
$
and 
$
\sum_{\tau=1}^{\infty}f_k(u(\tau))=\sum_{\tau=A_k+1}^{A_k+W_k}f_k(u(\tau)).  
$ 
 Similar to  Theorem~\ref{thm:LL} one can show that $H=\lambda G$ applies to all {\em coherent} systems.

\begin{cor}
Under condition L2, $H=\lambda G$ applies to all {\em coherent} systems.
\end{cor}
\begin{cor} Relation (\ref{eq:fktau}) applies to all {\em coherent} systems. In particular,
the results in Subsection~\ref{subsection:ggc} apply to all {\em coherent} queues.
\end{cor}
{\bf Proof.}
 Let 
\bear\label{eq:fkutau}
f_k(u(\tau)) &=& S_k {\bf 1}\{ A'_k < u(\tau) \leq (A_k+W_k^q)' \} +\nonumber                 \\
       & &  (S_k-(\tau-A_k-W_k^q))
               {\bf 1}\{ (A_k+W_k^q)' < u(\tau) \leq (A_k+W_k^q+S_k)' \} \; .
\end{eqnarray}
The result follows by
noting that for all {\em coherent} systems we have,
\[
\sum_{\tau=1}^{\infty}f_k(u(\tau))=\sum_{\tau=1}^{\infty}f_k(\tau)=G_k \;;
\]
where $f_k(\tau)$ and $f_k(u(\tau))$ are given by (\ref{eq:fktau}) and (\ref{eq:fkutau}) respectively.\done
\bigskip

\noindent
{\bf Remarks.}

(i) Consider the $B/Geom/1$ {\em coherent} model.  Assume arrival probability is $\alpha$ and mean service time is $1/\beta$. Then by (\ref{eq:spPK}), we have
\[
W_q=\rho(1-\beta)/(\beta-\alpha)\;,
\]
and
\[
W=(1-\alpha)/(\beta-\alpha)
\]

(ii) The workload process goes up by $S_k, k=1,2,\ldots ,$ at arrival instants, then it goes down by one unit at a time.  This this equivalent to  a model with batch arrivals where the batch size is equal to $S$ and service times are deterministic with  one unit per customer. Then the number of customers in the system at any given time is equal to corresponding workload $V(\tau)$. Therefore, the mean  queue-length is the same as the mean work in the system  and can be computed from the  $PK$ formula given in this section. 

(iii) The $PK$ formula given in (\ref{eq:spPK})  has been derived  for one {\em SR} rule at a time. See for instance  Hunter~\cite{Hun83b} who gives the result for {\em LAS-DA} observed at departure instants, and Chan et al.~\cite{Cha97} who gives the same result for discretized $M/G/1$  queue using the equivalent of the  {\em EAS} scheduling rule. In contrast, our approach applies to all {\em coherent} systems.

\section{Concluding Remarks}\label{sec:cr}

The results in this article open the door to study discrete-time queues as classes.
By studying one model within  a class (e.g., {\em coherent} systems), one can draw conclusions about the entire class. More research is needed to study the invariant characteristics  of each class, with emphasis on the {\em coherent} class. 
When using {\em SR}, one also needs to be aware of the corresponding observation epoch.  
Although most of  systems studied in the literature are {\em coherent}, there  are instances when {\em incoherent} systems can be useful as they can give insights about the behavior of a {\em coherent} system. An example is when we apply {\em BASTA} (El-Taha~\cite{Elt24}) for discrete-time systems with {\em SR}.

\section{Appendix: Proof of  Little's Law}

We shall need the following lemma, which is a  discrete-time sample-path analogue of the elementary renewal theorem, was proved in  El-Taha and Stidham \cite{Elt99} 
in continuous time.
\begin{lem} \label{lem:nt}
Let $0 \leq \lambda \leq \infty$. 
Then $\tau^{-1}A(\tau) \rightarrow \lambda$ as $\tau\rightarrow \infty$ 
if and only if $k^{-1}A_k \rightarrow \lambda^{-1}$ as
$k \rightarrow \infty$.
\end{lem}
\noindent
For a proof of Lemma~\ref{lem:nt} refer to the proof Lemma 2.1 of El-Taha and Stidham \cite{Elt99}.
Now, we obtain the basic inequality,
\begin{lem}\label{lem:ineq} For all $\tau \geq 1$
\beq   \label{eq:bbasic}
\sum_{k:A_k \leq \tau} W_k \geq \sum_{j=1}^\tau L(j) 
                        \geq \sum_{k:D_k \leq \tau} W_k \; , \; \tau \geq 1 \; .
\eeq
\end{lem}
{\bf Proof.}
Using the above notation, we obtain the  basic equalities,
\begin{eqnarray}
\sum_{j=1}^\tau L(j) &=& \sum_{k:A_k \leq \tau} W_k - 
\sum_{k:A_k < \tau \leq D_k} (D_k-\tau)\;;\\ 
\sum_{j=1}^\tau L(j) &=& \sum_{k:D_k \leq \tau} W_k + 
\sum_{k:A_k < \tau \leq D_k} (\tau-A_k)\;;\\
\sum_{k:A_k < \tau \leq D_k} W_k &=& 
	\sum_{k:A_k \leq \tau} W_k - \sum_{k:D_k \leq \tau} W_k \hbox{ .} 
\end{eqnarray} 
Note that the third equality follows from the first two. The Lemma then follows.\done

\noindent
{\bf Proof of Theorem~\ref{thm:lastword}.}
 Since  $W < \infty$,
it follows that $n^{-1}W_n \rightarrow 0$. 
 Moreover,  $W_n/A_n \rightarrow 0$ as $n \rightarrow \infty$. To see this write $W_n/A_n=(W_n/n)(n/A_n)$ and use Lemma~\ref{lem:nt}.

  Let $\epsilon > 0$ be given. Since $W_n/A_n \rightarrow 0$ as $n \rightarrow \infty$,  there exists an integer $N$ such that, $k \geq N$
implies $W_k \leq A_k \epsilon$ .  Therefore, for all $\tau \geq 0$ ,
\begin{eqnarray*}
\sum_{k:D_k \leq \tau} W_k &=& \sum_{k:A_k+W_k \leq \tau} W_k \\
                        &\geq& \sum_{k \geq N: A_k(1+\epsilon) \leq \tau} W_k \\
                        &\geq& \sum_{k:A_k(1+\epsilon) \leq \tau} W_k -
                                     \sum_{k \leq N-1} W_k ,
\end{eqnarray*}
which, together with the basic inequality~(\ref{eq:bbasic}), implies
\begin{eqnarray}
\sum_{k:A_k \leq \tau} W_k \geq \sum_{j=1}^{\tau} L(j) \label{eq:AgeqL}% \\
                        \geq \sum_{k:A_k(1+\epsilon) \leq \tau} W_k -
                                     \sum_{k \leq N-1} W_k \;. \label{eq:DgeqA} 
\end{eqnarray}
Moreover,
\begin{eqnarray*}
\limtau \tau^{-1} \sum_{k:A_k \leq \tau} W_k
         = \limtau \tau^{-1} A(\tau)A(\tau)^{-1} \sum_{k=1}^{A(\tau)} W_k 
         =\lambda W \;.
\end{eqnarray*}

and 
\begin{eqnarray*}
\limtau \tau^{-1} \sum_{k:A_k(1+\epsilon) \leq \tau} W_k 
    &=& (1 + \epsilon)^{-1} \limtau [\tau(1 + \epsilon)^{-1}]^{-1} 
             \sum_{k:A_k \leq \tau(1+\epsilon)^{-1}} W_k \\
    &=& (1 + \epsilon)^{-1} \limtau \tau^{-1} \sum_{k:A_k \leq \tau} W_k \\
    &=& (1 + \epsilon)^{-1} \lambda W\;.
\end{eqnarray*}
Now,
\begin{eqnarray*}
\lambda W=\lim_{\tau \rightarrow \infty} \tau^{-1} \sum_{k:A_k \leq \tau} W_k
    &\geq&  \limsup_{\tau\rightarrow\infty} \tau^{-1} \sum_{j=1}^{\tau} L(j) \\
&\geq&  \liminf_{\tau\rightarrow\infty} \tau^{-1} \sum_{j=1}^{\tau} L(j) \\
    &\geq& \lim_{\tau \rightarrow \infty} \tau^{-1} \sum_{k:A_k(1+\epsilon)
                 \leq \tau} W_k - \limtau\tau^{-1} \sum_{k \leq N-1} W_k\\
    &=& (1 + \epsilon)^{-1} \lambda W\;.
\end{eqnarray*}
 Since $\epsilon > 0$ was arbitrary,  using the fact that $\limtau \tau^{-1} \sum_{k \leq N-1} W_k = 0$, we conclude  that these inequalities hold in the limit as $\epsilon \rightarrow
0$. Therefore, the limit $\limtau \tau^{-1} \sum_{j=1}^{\tau} L(j)$ exists, and 
\[
\lambda W\geq  \limtau \tau^{-1} \sum_{j=1}^{\tau} L(j)\geq \lambda W\;.
\]
This completes the proof. \done 
%%%%%%%%%%%%%%%%%%%%%%%%%%%%

%\bibliography{ref}

%\end{document}

%\bibliography{ref}

\end{document}